\PassOptionsToPackage{numbers,sort&compress}{natbib}
\documentclass[smallextended,natbib]{svjour3}

\usepackage{graphicx}
\usepackage{amsmath}
\usepackage{txfonts}
\usepackage{siunitx}
\usepackage{booktabs}
\usepackage{multirow}
\usepackage{xcolor}
\usepackage{soul}
\usepackage{nicefrac}
\usepackage{bbold}
\usepackage{tikz}
\usetikzlibrary {math} 
\usepackage{tikz-3dplot}
\usepackage[hidelinks]{hyperref}
\usepackage{xurl}

\ifdefined\qty\else
  \ifdefined\NewCommandCopy
    \NewCommandCopy\qty\SI
  \else
    \NewDocumentCommand\qty{O{}mm}{\SI[#1]{#2}{#3}}
  \fi
\fi

\smartqed
\setcitestyle{numbers,square,comma}
\journalname{Journal of Scientific Computing}

\makeatletter\def\makeheadbox{}\makeatother

\usepackage[a4paper,left=3.2cm,right=3.2cm,top=2.6cm,bottom=2.8cm]{geometry}

\begin{document}

\title{Variational Inference Using a Differentiable Multigrid Linear Solver}
\titlerunning{VI Using a Differentiable Multigrid Linear Solver}

\author{%
  Andr\'es~Ram\'irez \and
  Philipp~Haim \and
  David~F\"oger \and
  Torsten~A.~En{\ss}lin \and
  Philipp~Frank \and
  Philipp~Gschwandtner \and
  Dominik~J\"ustel \and
  Philipp~Mertsch \and
  Vo~Hong~Minh~Phan \and
  Laurin~S\"oding \and
  Hanieh~Zandinejad \and
  Ralf~Kissmann%
}
\authorrunning{A. Ram\'irez et al.}

\institute{%
  A.~Ram\'irez, R.~Kissmann \at
  Universit\"at Innsbruck, Institut f\"ur Astro- und Teilchenphysik,
  Technikerstr.~25/8, 6020 Innsbruck, Austria \\
  \email{Andres.Ramirez-Tapias@student.uibk.ac.at}
  \and
  P.~Haim, D.~J\"ustel \at
  Institute of Biological and Medical Imaging, Institute of Computational
  Biology, and Institute of AI for Health, Bioengineering Center and
  Computational Health Center, Helmholtz Zentrum M\"unchen,
  D-85764 Neuherberg, Germany \\
  Chair of Biological Imaging, Central Institute for Translational Cancer
  Research (TranslaTUM), School of Medicine and Health \& School of
  Computation, Information and Technology, Technical University of Munich,
  Munich, Germany
  \and
  D.~F\"oger, P.~Gschwandtner \at
  Universit\"at Innsbruck, Institut f\"ur Informatik,
  Technikerstr.~21\,a, 6020 Innsbruck, Austria
  \and
  H.~Zandinejad \at
  Max Planck Institute for Astrophysics,
  Karl-Schwarzschild-Stra{\ss}e~1, 85748 Garching bei M\"unchen, Germany \\
  Ludwig Maximilian University of Munich,
  Geschwister-Scholl-Platz~1, 80539 M\"unchen, Germany
  \and
  T.~A.~En{\ss}lin \at
  Max Planck Institute for Astrophysics,
  Karl-Schwarzschild-Stra{\ss}e~1, 85748 Garching bei M\"unchen, Germany \\
  Ludwig Maximilian University of Munich,
  Geschwister-Scholl-Platz~1, 80539 M\"unchen, Germany \\
  Deutsches Zentrum f\"ur Astrophysik, Postplatz~1, 02826 G\"orlitz, Germany \\
  Excellence cluster ORIGINS, Boltzmannstr.~2, 85748 Garching, Germany
  \and
  P.~Frank \at
  Max Planck Institute for Astrophysics,
  Karl-Schwarzschild-Stra{\ss}e~1, 85748 Garching bei M\"unchen, Germany
  \and
  P.~Mertsch \at
  Institute for Theoretical Particle Physics and Cosmology,
  RWTH Aachen University, Sommerfeldstr.~16, 52074 Aachen, Germany
  \and
  L.~S\"oding \at
  Institute for Theoretical Particle Physics and Cosmology,
  RWTH Aachen University, Sommerfeldstr.~16, 52074 Aachen, Germany \\
  Max Planck Institute for Astrophysics,
  Karl-Schwarzschild-Stra{\ss}e~1, 85748 Garching bei M\"unchen, Germany
  \and
  V.~H.~M.~Phan \at
  Sorbonne Universit\'e, Observatoire de Paris, PSL Research University,
  LERMA, CNRS UMR~8112, 75005 Paris, France
}

\date{}

\maketitle

\begin{abstract}
Gradient-based Bayesian inference methods require efficient access to Jacobian and adjoint-Jacobian operators of high-dimensional forward models. While multigrid solvers provide near-optimal complexity for elliptic partial differential equations, they are rarely available in forms compatible with automatic differentiation (AD).
We develop a differentiable multigrid solver for steady-state diffusion--absorption problems and derive its adjoint operations analytically through the full multigrid hierarchy. The resulting solver, DMGS, is implemented in C++ and interfaced with JAX to provide efficient Jacobian--vector and vector--Jacobian products for variational inference in the NIFTy framework.
We validate the approach on a 3D inverse problem involving diffuse radiative transfer in tissue, reconstructing an effective radiative source from Monte Carlo-simulated data. The reconstruction reproduces the data at a reduced chi-squared of $1.1$ and generalizes to 32 independent validation datasets.
Benchmarks against a JAX-native multigrid implementation show comparable runtimes and consistently lower peak memory for the hand-derived adjoint, with modest reverse-mode overhead.
These results establish differentiable multigrid solvers as practical building blocks for variational inference in PDE-constrained problems.
\keywords{Multigrid methods \and Automatic differentiation \and Adjoint methods \and Variational inference \and Partial differential equations}
\end{abstract}

\section{Introduction}
\label{Introduction}

Many scientific computing applications couple large-scale PDE solvers with gradient-based Bayesian inference, which requires efficient and accurate derivatives in high-dimensional parameter spaces \citep{Goldman02102022,THIAGARAJAN2025118401}. When the forward model is defined by partial differential equations (PDEs), discretization leads to large linear systems, and computing Jacobians and adjoint-Jacobian products can dominate the total computational cost of inference \citep{griewank2008evaluating, baydin2018automatic}. 

Automatic differentiation (AD) provides machine-precision derivatives at a computational cost comparable to evaluating the forward model itself. Consequently, AD has become a core technology in probabilistic programming languages (PPLs) such as Stan \citep{stan2017}, PyMC3 \citep{pymc2023}, Pyro \citep{bingham2019pyro}, NumPyro \citep{phan2019numpyro}, and Edward \citep{tran2016edward}, which have emerged as AD-ready toolkits for accessible and efficient derivative evaluations in statistical modeling. These frameworks integrate gradient-based inference algorithms using forward- and reverse-mode AD. Such algorithms include Automatic Differentiation Variational Inference (ADVI) \citep{kucukelbir2017advi} and Stochastic Variational Inference (SVI) \citep{svi2013}. In inverse problems and variational inference (VI) schemes, these operators are important for expressing gradient updates. Reverse-mode AD, in particular, provides efficient access to adjoint-Jacobian products $J^\dagger w$, where $J$ is the Jacobian of the forward model and $w$ a cotangent vector in data space, at a computational cost that is essentially independent of the number of input parameters. 

Iterative multigrid (MG) methods are widely used in state-of-the-art solvers for linear elliptic and parabolic PDEs, and offer near-optimal $O(N)$ complexity in the number of unknowns $N$ of the discretized system by combining relaxation schemes with hierarchical grid corrections \citep{briggs2000multigrid, trottenberg2001multigrid}.
For PDE-constrained problems, developments such as the adjoint-state method \citep{giles2000adjoint} provide a general and computationally efficient route to obtain gradients, with a cost essentially equal to one forward plus one reverse computation. However, such methods are often defined for performance and are not written in forms compatible with inference frameworks (i.e., adjoint-Jacobian products are not readily available). Furthermore, MG solvers rely on recursively defined inter-grid transfers, coarse-grid corrections, and specialized smoothers that are typically tuned for performance and are not expressed in differentiable programming paradigms. Recent developments in differentiable programming have begun to bridge this gap. Automated adjoint frameworks for finite-element and PDE problems, such as dolfin-adjoint \citep{farrell2013dolfinadjoint}, demonstrate how to derive Jacobian and adjoint-Jacobian models automatically from high-level problem specifications. Furthermore, AD-native libraries such as JAX \citep{jax2018github} provide composable transformations (e.g., vector-Jacobian and Jacobian-vector products) that make it straightforward to implement differentiable solvers, including multigrid, in Python.

Recent developments in JAX and its ecosystem of transformations make it straightforward to compose solver primitives (linear solvers, multigrid cycles, and time integrators) into fully differentiable models that expose both Jacobian--vector products $Jv$, with $v$ a tangent vector in parameter space, and adjoint-Jacobian products $J^\dagger w$ \citep{yashchuk2023bringing}. Continued efforts have been made to implement AD in JAX-based simulators: BlackBIRDS \citep{Quera-Bofarull2023} provides automatic inference primitives for differentiable simulators implemented in PyTorch; MadJax \citep{heinrich2022madjax} aims to produce differentiable components (e.g., matrix elements in specific physics domains) for gradient-based learning; other packages, such as JAX-DIPS \citep{mistani2022jaxdips} and JAX-Fluids \citep{bezgin2023jaxfluids}, demonstrate high-dimensional differentiable simulations useful for inverse problems and learning tasks. Moreover, there have been multiple developments in AD-compatible simulators outside of the JAX ecosystem as well, although these efforts remain focused on optimization problems \citep{newbury2024review}. Despite the versatility of automatic derivatives, there are trade-offs: AD-native implementations enable rapid software prototyping and straightforward correctness thanks to automatic adjoint operators, whereas hand-coded native adjoint operators can offer substantially better raw performance and memory behavior for very large-scale runs.

Despite the maturity of multigrid methods and AD frameworks individually, their integration within VI remains underexplored, particularly in uniting efficient differentiable solvers with structured VI schemes that go beyond mean-field approximations. Existing work provides either AD-native multigrid prototypes or hand-optimized forward solvers, but rarely a solver with explicit adjoint operators that is directly usable inside metric-based VI pipelines.

The main contribution of this work is therefore to show that PDE multigrid adjoints are practically available for VI: we derive an adjoint multigrid solver (DMGS) for diffusion--absorption problems through the full hierarchy; expose it as a JAX-compatible differentiable building block; and integrate it into the NIFTy framework \citep{nifty1,nifty3,nifty5,niftyre} to perform metric-based VI on a 3D reconstruction problem, using Metric Gaussian Variational Inference \citep[MGVI,][]{knollmueller2019mgvi} and geometric Variational Inference \citep[geoVI,][]{frank2021geovi}. A major technical challenge is to differentiate consistently through the recursive transfers, coarse-grid corrections, and smoothers so that the Jacobian and adjoint-Jacobian operators required by these schemes are well defined. Both MGVI and geoVI have been demonstrated to provide accurate results in complex field reconstruction problems, as well as to quantify the remaining uncertainty in terms of (approximate) posterior samples, with recent applications in astrophysics \citep{tsouros2024,gasreco2025,edenhofer2024} and other domains \citep{2023LimOM21734R,2025SciA117319B,2025arXiv251202204H}.

We demonstrate the approach on a medical imaging problem involving diffuse radiative transfer in tissue, where the goal is to reconstruct an effective radiative source from synthetic observations, as measured, for example, with optoacoustic devices. We focus on a multigrid solver (MGS) for diffusion--absorption problems, particularly the same implementation as introduced for the problem of Galactic Cosmic-Ray Transport with the PICARD solver \citep{KISSMANN201437}. We exploit its mathematical properties to define the operations $S$ and $S^\dagger$ that define DMGS and utilize them in the context of VI. To quantify the cost of hand-derived adjoints, we benchmark DMGS, implemented in C++, against a JAX-native implementation relying on automatic differentiation (JMGS). Details regarding JMGS will be presented in a future publication; here its role is limited to this comparison. To our knowledge, this work is among the first demonstrations of differentiable multigrid driving metric-based VI.

The remainder of this work is structured as follows: in Section~\ref{sec:Methodology}, we introduce all the ingredients for the construction of DMGS and its corresponding adjoint implementation, as well as the specifics on the VI setup; in Section~\ref{Results}, we analyze the results of our reconstruction and compare the performance of the forward- and reverse-modes of DMGS and JMGS; in Section~\ref{Conclusions} we present our conclusions and outlook.

\section{Methodology}
\label{sec:Methodology}

In this section we describe the construction of the differentiable multigrid solver and its adjoint operations. We first introduce the diffusion--absorption forward model, then derive the multigrid solver and its adjoint operators, wrap DMGS for use in NIFTy, and finally specify the radiative-transfer reconstruction setup. We begin with the PDE (see Section~\ref{subsec:Jacobian} and Fig.~\ref{fig:MGSgraph}). Let $\mathcal{L}$ be a linear differential operator acting on a scalar field $\rho$, and consider the steady-state diffusion--absorption equation

\begin{equation}
\label{maineq}
    \mathcal{L}\rho=\left(-\nabla \cdot(D\cdot\nabla)+\lambda\right) \rho = q.
\end{equation}

Here, $\rho$ and $q$ denote the fluence and source fields, respectively; $D$ is a spatially varying isotropic diffusion coefficient (equivalently, a diagonal diffusion tensor $D\,\mathbb{1}$, with $\mathbb{1}$ the identity); and $\lambda$ is the absorption coefficient.

After discretization on a structured grid, $\mathcal{L}$ becomes a linear operator $L$ acting on coefficient vectors in $\mathbb{R}^{n_{\mathrm{dof}}}$. The explicit form of $L$ follows from second-order finite differences as introduced in \citep{KISSMANN201437}. Solving the system corresponds to applying the inverse operator $L^{-1}(q)=\rho$. In practice, this inverse is approximated numerically by a solver $S$, which defines the mapping

\begin{equation}
    \rho = S(q)
\end{equation}
and is represented by a composition of linear operators (see Section~\ref{subsec:SolverDef} for details). To implement the solver $S$ in a probabilistic programming language (PPL), we need to define a push-forward operation and a pull-back operation (respectively, $S$ and $S^\dagger$; see Section~\ref{subsec:Jacobian} for details).

With a PPL-compatible solver $S$, we study its applicability in the context of inference. For this purpose, we use the Gaussian process (GP) models and VI algorithms implemented in the Bayesian imaging library NIFTy, particularly the JAX implementation NIFTy.re. To enable the use of DMGS within NIFTy, we wrap this solver in JAX using the library JAXbind \citep{Roth2024}.

\subsection{VI with NIFTy}
\label{subsec:VINIFTy}

In NIFTy, VI algorithms approximate the posterior using an approximate distribution function, such as a multivariate Gaussian, which allows one to draw approximate posterior samples. This is done by minimizing the Kullback--Leibler (KL) divergence between the variational approximation and the exact posterior distribution. The core ingredients of these approximations are: the likelihood $\mathcal{P}(d|\theta)$ of a dataset $d$ given the standardized parameters $\theta$, a forward model $f(\theta)$ mapping $\theta$ to the data space, and the Jacobian

\begin{equation}
    J(\theta)=\frac{df}{d\theta}
\end{equation} 
of the mapping. Both MGVI and geoVI rely on $J(\theta)$ to define the local information geometry of the posterior through the standardized Fisher metric

\begin{equation}
    I(\theta) = J(\theta)^\dagger I_d(\theta) J(\theta) + \mathbb{1}.
\end{equation}

Here, $I_d(\theta)$ is the Fisher metric of the model likelihood. The inverse of the standardized Fisher metric provides a curvature approximation in latent space and determines the local covariance of the variational distribution, $\Theta(\theta) \approx I(\theta)^{-1}$. All core steps involving the Fisher metric and coordinate transformations of the parameters (i.e., sampling, standardizing, updating by gradient-based updates) in MGVI and geoVI require repeated evaluations of the Jacobian--vector (push-forward) and vector--Jacobian (pull-back) products. NIFTy.re performs these evaluations through implicit operators via the corresponding JAX \texttt{jvp} and \texttt{vjp} calls. Therefore, a forward model must possess a JAX-compatible definition of $J(\theta)$ and $J(\theta)^\dagger$, both of which are computed via the chain rule for complex models composed of multiple mappings. In the context of this work, one of these mappings is performed by a multigrid solver. Thus, we need to define the Jacobian and Jacobian-transpose of such a solver, using algorithmic differentiation.

\subsection{Structure of the solver}
\label{subsec:SolverDef}

The multigrid solver $S$ studied here can be represented as a composition of $l$ linear operations $g_i$ in the form of

\begin{equation}
\label{Eqcomp}
    S = g_l\circ g_{l-1}\circ\cdots\circ g_1.
\end{equation}

In our notation, $S_j$ where $j\leq l$ is the $j$-th composition $g_j\circ\cdots\circ g_1$, and $S=S_l$.

The solver operates on a hierarchy of $n_{\mathrm{lev}}$ Cartesian grids. Level 0 corresponds to the finest grid, while level $n_{\mathrm{lev}}-1$ represents the coarsest grid. The grid resolution is reduced by a factor of two in each dimension between consecutive levels, such that 

\begin{equation}
    h^i_k=2h^i_{k-1},
\end{equation}
with $h^i_{k}$ the respective grid-point distance at a level $k$ along the dimension $i$. In this work we follow the notation from \citep{trottenberg2001multigrid} and define our multigrid framework as follows: We assume equal spacing in all dimensions and choose the number of grid points per dimension at the finest level to be $2^{n_g}+1$, which determines $n_{\mathrm{lev}}$. Each coarse grid is constructed by halving the resolution of the preceding level. We derive the discrete form $L^k$ of Eq.~(\ref{maineq}) using finite differences for each level, where the solution at level 0 is calculated by a multigrid cycle (i.e., V--cycle for our application), using the Red--Black Gauss--Seidel (RB-GS) relaxation method and assuming Dirichlet boundary conditions. In our setup, the cycle is repeated 15 times to ensure a solution with sufficiently high precision.

\begin{figure}
    \centering
    \includegraphics[width=\linewidth]{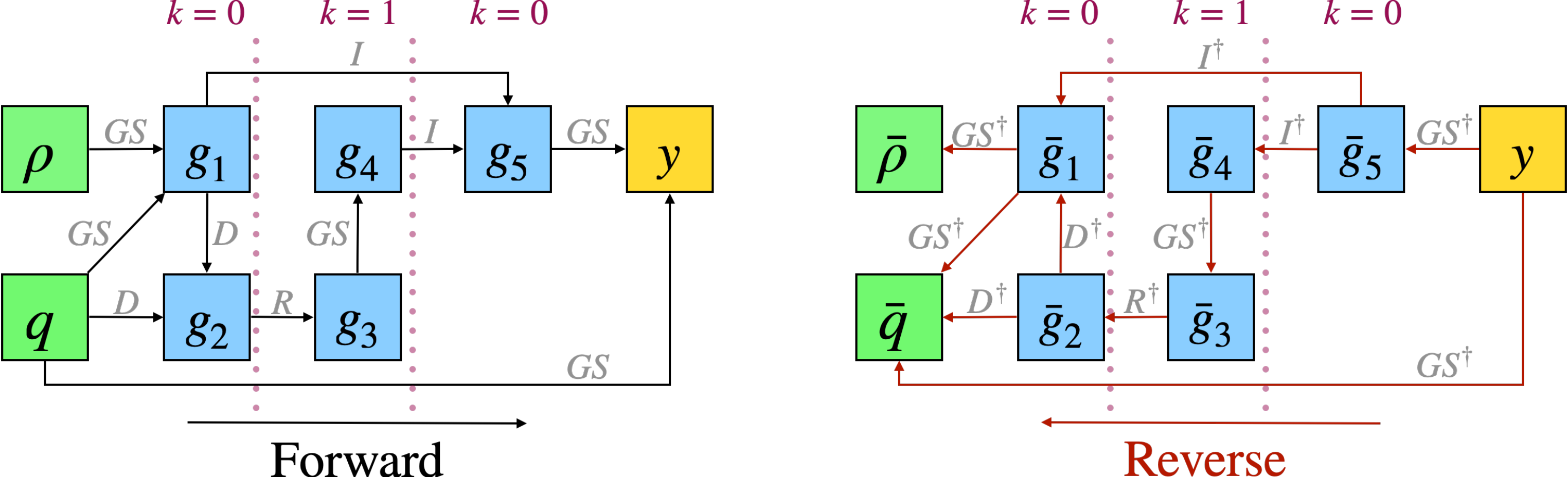}
    \caption{Computational graph for a 2-level RB-GS multigrid solver for its forward (left) and reverse (right) modes, connecting an initial state $\rho$ and a source term $q$ to the final solution of $S$. Here, the $g_i$ and $\bar{g}_i$ nodes represent the corresponding output of the linear operations for the forward and reverse modes, respectively, where $i$ denotes the order in which the operations are performed in the forward case. The operations performed at each level $k$ are grouped by dotted lines. The upper row of nodes corresponds to the intermediate states of all $\rho_k$, while the lower row of nodes corresponds to $q_k$. The arrows connecting all nodes are labeled based on the operations being performed: $GS$ for the Red--Black Gauss--Seidel relaxation, $I$ for interpolation (prolongation), $D$ for the defect, and $R$ for the restriction. The symbol $\dagger$ denotes the adjoint implementation.}
    \label{fig:MGSgraph}
\end{figure}

We illustrate the series of operations that compose our multigrid implementation for a 2-level example in Fig.~\ref{fig:MGSgraph}. There are four major operations taking place in the solver: an RB-GS operation (see Section~\ref{subsec:RBGS}); a defect operation (see Section~\ref{subsec:defect}); and two transfer operations (see Section~\ref{subsec:Transfer}): prolongation (onto a finer grid) and restriction (onto a coarser grid). The RB-GS operation maps, e.g., $GS(\rho,q)=g_1$, where the even (red) grid points of $g_1$ are updated first, and the odd (black) grid points follow. The defect operation, e.g., $D(g_1,q)=g_2$, defines the error between $g_1$ and $q$, which is used as a new source term after it is transferred to the next grid level via the restriction operation. For prolongation and restriction, we use full-weighting operators; restriction is defined by the following stencil

\begin{equation}
\label{eq:rest}
R_{h\rightarrow H} = \frac{1}{4}
\begin{bmatrix}
1 & 2 & 1
\end{bmatrix}
\end{equation}

\subsection{Defining the Jacobian}
\label{subsec:Jacobian}

To construct the necessary push-forward and pull-back operations in the context of a differentiable solver, we compute the gradients of the operator $S$ with respect to the source density $q$ and follow the chain rule using the notation of Eq.~(\ref{Eqcomp})

\begin{equation}
    \frac{\partial S_j}{\partial q}=\sum_{i\prec j}g_j'\frac{\partial g_i}{\partial q},
\end{equation}
where $g_j'=\partial g_j/\partial g_{i}$. The notation $i \prec j$ denotes all nodes $i$ which are predecessors of the node $j$ in Fig.~\ref{fig:MGSgraph}. Following the graph, to construct the Jacobian of $S$, it is sufficient to compute the partial derivatives of every intermediate function $g_j$ with respect to its corresponding predecessors $g_i$. Furthermore, because $S$ is linear in $q$ for fixed optical coefficients $D$ and $\lambda$, fixed Dirichlet boundary conditions, and a fixed multigrid schedule, the Jacobian of $S$ equals $S$ itself: a push-forward Jacobian--vector product applies $S$ directly to the perturbation. For the pull-back operation, the same derivatives $g_i'$ need to be calculated. However, the $i$-th pull-back operation $S^\dagger_i$ is applied starting from the successors of $S_i$, instead of the predecessors. Therefore, the order in which the operations $g_j$ are applied is reversed. In contrast to $S$, the definition of $S^\dagger$ is not trivial in this case. Thus, we will focus on describing $S^\dagger$ in the context of our RB-GS MGS. For illustration purposes, we define the corresponding pull-back operations of the relevant $g_j$ components in one dimension only. For higher dimensions and off-diagonal terms in Eq.~(\ref{maineq}), the operations include all combinations of the corresponding neighbors in the other dimensions.

\subsubsection{RB-GS}
\label{subsec:RBGS}

The RB-GS relaxation method is a highly parallelizable variant of the Gauss--Seidel relaxation method for solving elliptic PDEs on structured grids. In one dimension, an application of the operation $GS(\rho^0_k,q_k)=\rho_k$, with $\rho^0_k$ the initial state and $\rho_k$ an updated state of the final solution $\rho$, has the form

\begin{align}
    \rho_k(x_i) &= \left[L^k_{i-1}\rho^0_k(x_{i-1})+L^k_{i+1}\rho^0_k(x_{i+1})-q_k(x_i)\right]/L^k_i
\end{align}

Here, $L^k_i$ is the discretization of the operator $L$ on level $k$ at position $x_i$. We note that the Red--Black convention means that odd grid points are calculated using the updated state of the even ones.
The adjoint operation $GS^\dagger(\rho_k)=(\bar{\rho}_k,\bar{q}_k)$ has the following form:

\begin{align}
    \bar{\rho}_k(x_i) &= L^k_i\left[\rho_k(x_{i-1})/L^k_{i-1}+\rho_k(x_{i+1})/L^k_{i+1}\right], \\
    \bar{q}_k(x_i) &= -\rho_k(x_i)/L^k_i.
\end{align}

On coarse levels ($k\neq0$), $\rho_k$ is initialized to zero. In the reverse-mode graph (Fig.~\ref{fig:MGSgraph}), $GS^\dagger$ therefore updates only one node per level, and the cotangent $\bar{\rho}_k$ contributes to the pull-back of $\bar{\rho}$ and $\bar{q}$ exclusively at $k=0$.

\subsubsection{Defect}
\label{subsec:defect}

For levels $k\neq0$, the solution $\rho_{k}$ calculated via the RB-GS operation measures the error of the expected solution given the source term $q_{k}$, also known as the defect $e_k$, for the discretized PDE in Eq.\,(\ref{eq:error}). 

\begin{equation}
\label{eq:error}
    L^k \rho_k  =  q_{k}
\end{equation}

The defect operation then calculates the remaining error $e_k$ between the solution $\rho_k$ and $q_k$, which is transferred to the coarser level $k+1$ via the restriction operator, defining $q_{k+1}$ as shown in Eq.\,(\ref{eq:defect}).

\begin{equation}
    \label{eq:defect}
    q_k-L^k\rho_k=e_k \rightarrow R(e_k) = q_{k+1}
\end{equation}

This way, the RB-GS operation rapidly corrects for the error at the scales defined by the different grid coarsenings via $k$. In a multigrid setup, the sequential application of the defect operation across all levels allows for the efficient smoothing of errors of the PDE on all scales. In one dimension, the defect operation $D(\rho_k,q_k)=e_k$ defines the error explicitly as

\begin{equation}
    e_k(x_i) = L^k_i \rho_k (x_i) - L^k_{i-1}\rho_k(x_{i-1})-L^k_{i+1}\rho_k(x_{i+1})+q_k(x_i).
\end{equation}

For the adjoint operation $D^\dagger(e_k)=(\bar{\rho}_k,\bar{q}_k)$, this means

\begin{align}
    \bar{\rho}_k(x_i) &= L^k_{i}\left[e_k(x_i)-e_k(x_{i-1})-e_k(x_{i+1})\right] \\
    \bar{q}_k(x_i) &= e_k(x_i)
\end{align}

\subsubsection{Transfer Operations}
\label{subsec:Transfer}

In a multigrid solver, the transfer operations define how relevant quantities (i.e., $\rho_k$ and $q_k$) are sampled on grids at different levels. For regular grids, the prolongation (from coarse to fine) and restriction (from fine to coarse) transfer operations between them are related via

\begin{equation}
    I_{H\rightarrow h}=2^nR_{h\rightarrow H}^T
\end{equation}

Here, $n$ is the dimensionality of the problem. This property holds when the restriction $R$ is a full-weighting operator as defined in Eq.\,(\ref{eq:rest}), defining the prolongation (interpolation) $I$. Therefore, we can conveniently write the pull-back operations $I^\dagger(\rho_k)=2^nR(\rho_k)=\bar{\rho}_{k+1}$ and $R^\dagger(q_{k+1})=2^{-n}I(q_{k+1})=\bar{\delta}_k$.

\subsection{Full adjoint operator}
\label{subsec:fulladjoint}
Combining the adjoint forms of all solver components yields the adjoint multigrid operator $S^\dagger$, obtained by reversing the order of operations. We now proceed to explicitly define the algebraic adjoint $S^\dagger$. Given that the operators $g_j$ are linear mappings for all $j$, we find from Eq.~(\ref{Eqcomp})

\begin{equation}
\label{EqcompT}
    S^\dagger=\bar{g}_1 \circ \bar{g}_2\circ\cdots\circ \bar{g}_l.
\end{equation}

Reversing the order of operations, as shown in Fig.~\ref{fig:MGSgraph} and Eq.~(\ref{EqcompT}), we represent the full pull-back (adjoint) operation $S^\dagger(\bar{\rho})=\bar{q}$. Usually, in an MG approach, the push-forward application of this method iterates until the desired precision is obtained, and the total number of iterations that the solver performs is a priori unknown. Therefore, the pull-back operation can be performed in two modes. The first calculates a source term $q'$ by applying the discretized operator $L$ to $\rho$, which is flexible but subject to error at lower resolutions, and then solves to the desired precision while storing the total number of iterations. The second fixes the number of iterations for both forward and reverse modes, eliminating the need to calculate the source term. All experiments in this work use the second mode, with fifteen fixed V--cycles per solve, so that forward and adjoint calls are matched and JIT-compiled benchmarks remain deterministic.

\subsection{Reconstruction: Diffuse Radiative Transfer}

Optoacoustic medical imaging modalities measure the optical properties of tissue to gain insights into tissue composition and function \citep{englert2025revmod}. This setting provides a realistic test case for differentiable PDE solvers within variational inference. To faithfully recover those optical properties, an accurate model of light transport in tissue is needed. In the case of optoacoustic tomography, a diffusion--absorption equation, Eq.~(\ref{maineq}), can be used to approximate the radiative transfer equation \citep{Cheong1990,Tarvainen2005}, provided that the radiation scattering coefficient in the tissue is much larger than its radiation absorption coefficient $\lambda_a\ll\lambda_s$, that the boundaries of the solution grid are sufficiently far from the region of interest, and that the radiation source is sufficiently far from the tissue to ensure that diffusion dominates over the ballistic regime.

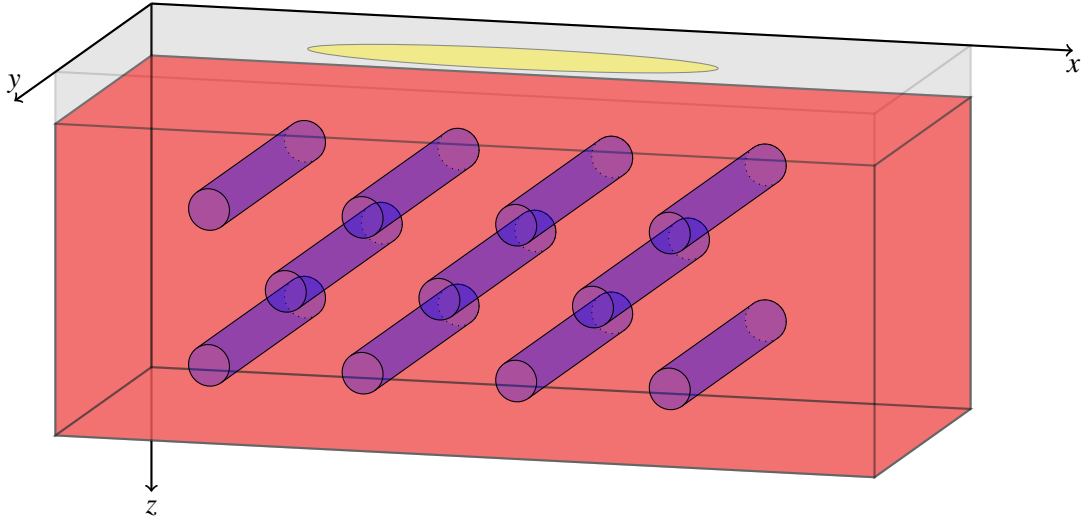
\begin{figure}
\centering
\tdplotsetmaincoords{79}{105}% theta,phi
\tdplotsetrotatedcoords{0}{-90}{0}
\begin{tikzpicture}[scale=1.4,tdplot_rotated_coords,
                    every node/.append style={font=\large}]
 
  \def\xmax{3.5}  
  \def\ymax{8}  
  \def\zmax{3.5} 
  \def\xinner{3.0} 
  
  \draw[thick] (3.5,0,3.5) -- (3,0,3.5);
  \draw[thick,opacity=0.4] (3.0,0,3.5) -- (-0.71,0,3.5);
  \draw[thick,->] (-0.71,0,3.5) -- (-1.2,0,3.5) node[anchor=north]{$z$};
  \draw[thick,->] (3.5,0,3.5) -- (3.5,9,3.5) node[anchor=north]{$x$};
  \draw[thick,->] (3.5,0,3.5) -- (3.5,0,-1.5) node[anchor=south]{$y$};

  % outer cube
  \draw[fill=gray,opacity=0.1,thick](0,0,0) -- (\xmax,0,0) -- (\xmax,\ymax,0) -- (0,\ymax,0) -- cycle; % xy lower
  \draw[fill=gray,opacity=0.1,thick](0,0,\zmax) -- (\xmax,0,\zmax) -- (\xmax,\ymax,\zmax) -- (0,\ymax,\zmax) -- cycle; % xy upper
  \draw[fill=gray,opacity=0.1,thick](0,0,0) -- (\xmax,0,0) -- (\xmax,0,\zmax) -- (0,0,\zmax) -- cycle; % xz lower
  \draw[fill=gray,opacity=0.1,thick](0,\ymax,0) -- (\xmax,\ymax,0) -- (\xmax,\ymax,\zmax) -- (0,\ymax,\zmax) -- cycle; % xz upper
  \draw[fill=gray,opacity=0.1,thick](0,0,0) -- (0,0,\zmax) -- (0,\ymax,\zmax) -- (0,\ymax,0) -- cycle; % yz lower
  \draw[fill=gray,opacity=0.1,thick](\xmax,0,0) -- (\xmax,0,\zmax) -- (\xmax,\ymax,\zmax) -- (\xmax,\ymax,0) -- cycle; % yz upper
  
  % inner cube
  \draw[fill=red,opacity=0.3,thick](0,0,0) -- (\xinner,0,0) -- (\xinner,\ymax,0) -- (0,\ymax,0) -- cycle; % xy lower
  \draw[fill=red,opacity=0.3,thick](0,0,\zmax) -- (\xinner,0,\zmax) -- (\xinner,\ymax,\zmax) -- (0,\ymax,\zmax) -- cycle; % xy upper
  \draw[fill=red,opacity=0.3,thick](0,0,0) -- (\xinner,0,0) -- (\xinner,0,\zmax) -- (0,0,\zmax) -- cycle; % xz lower
  \draw[fill=red,opacity=0.3,thick](0,\ymax,0) -- (\xinner,\ymax,0) -- (\xinner,\ymax,\zmax) -- (0,\ymax,\zmax) -- cycle; % xz upper
  \draw[fill=red,opacity=0.3,thick](0,0,0) -- (0,0,\zmax) -- (0,\ymax,\zmax) -- (0,\ymax,0) -- cycle; % yz lower
  \draw[fill=red,opacity=0.3,thick](\xinner,0,0) -- (\xinner,0,\zmax) -- (\xinner,\ymax,\zmax) -- (\xinner,\ymax,0) -- cycle; % yz upper

  % VARIABLES
  \def\L{3.5}    % detector length
  \def\R{0.2}   % detector cylinder radius

\foreach \ix in {1,...,3}{
	\def\iymax{4}
	\ifnum \ix=2
	\tikzmath{\iymax=3;}
	\fi

	\foreach \iy in {1,...,\iymax}{
		
		\tikzmath{\sycart=1.5*\iy;}
		\ifnum \ix=2
		\tikzmath{\sycart=\sycart+0.75;}
		
		\fi
		\tikzmath{\sxcart=0.75*\ix;}
  
	    \def\ang{54} % rotate lines to simulate cylinder

		\begin{scope}[rotate around z=\ang]
		\tikzmath{\sx=(\sxcart*cos(\ang)+\sycart*sin(\ang));\sy=(-\sxcart*sin(\ang)+\sycart*cos(\ang));}

		  \fill[color=blue, opacity=0.4] (\sx,\sy+\R,\L) --++ (0,0,-\L) arc(90:-90:\R) --++ (0,0,\L) arc(270:90:\R) -- cycle; % inner back of cylinder
		\draw[] (\sx,\sy,0)++(90:\R) --++ (0,0,\L); % top horizontal
  		\draw[] (\sx,\sy,0)++(-90:\R) --++ (0,0,\L); % bottom horizontal
		\tdplotdrawarc[fill=blue,fill opacity=0.3]{(\sx,\sy,0)}{\R}{0}{360}{}{} % transverse plane at z=L/2 (front rings)
  		\tdplotdrawarc[thin,dotted,fill=blue,fill opacity=0.3]{(\sx,\sy,\L)}{\R}{0}{360}{}{} % transverse plane at z=-L/2 (back ring)
  		\tdplotdrawarc[]{(\sx,\sy,\L)}{\R}{-90}{90}{}{} % transverse plane at z=-L/2 (outer back ring)
	\end{scope}

}
}

  \begin{scope}[canvas is yz plane at x=3.5,transform shape]
  \draw[fill=yellow, opacity=0.4] (4,1.75) ellipse (2 and 0.5);
  \end{scope} 

\end{tikzpicture}
\caption{\label{FigSketch}
Sketch of the medical-imaging application used in our analysis including our grid definition.
At the top there is a light source (indicated by the yellow ellipse), slightly above the illuminated tissue (in red).
Additionally, we consider some blood-vessel-like structures (indicated by the blue cylinders).}
\end{figure}

To validate DMGS within a Bayesian inference context, we consider diffuse light transport in tissue as a forward problem: given spatially varying optical coefficients and an effective radiative source, predict the spatial distribution of absorbed energy; we then infer the source from noisy observations using VI. In optoacoustic imaging, such source reconstruction is a prerequisite for recovering absorption maps from measured signals and requires a differentiable forward model with accurate adjoint operators \citep{haim2026deeplight}. As an illustration of this setting, Fig.~\ref{FigSketch} shows a light source illuminating tissue that contains blood-vessel-like tubular structures.
Here, the idea is that these structures have different material properties regarding the scattering of light, which allows one to identify them in a real-world application.
Working in the diffusion-dominated regime, this pertains in particular to differences in the diffusion coefficient and absorption coefficient. Under the diffusion approximation to the radiative transfer equation \citep{Cheong1990,Tarvainen2005}, the reduced scattering coefficient is $\lambda_s' = (1-g)\lambda_s$, with $g$ the scattering anisotropy. The isotropic diffusion coefficient and absorption coefficient entering Eq.~(\ref{maineq}) are then
\begin{equation}
\label{eq:medimcoeffs}
    D = \frac{1}{3(\lambda_a + \lambda_s')} \quad\text{and}\quad \lambda = \lambda_a,
\end{equation}
respectively; Fig.~\ref{fig:medimcoeffs} shows the resulting spatially varying scalar field $D(\vec{r})$.

Furthermore, the distribution $\rho$ corresponds to the radiative fluence within the tissue. Following \citep{englert2025revmod}, we define the absorbed energy density as $L_e=\lambda\rho$. Given known $\lambda_s$, $\lambda_a$, and $g$, and the fluence $\rho$ from the solver, we reconstruct the source term $q$ in Eq.~(\ref{maineq}).

Throughout this work, given the 3D nature of $\rho,q,D\text{ and }\lambda$, these quantities will be presented via 3 projections (i.e., orthogonal, depth, and lateral) along the different axes with respect to the $xz$-plane, namely the central image plane. Hence, the $y$-axis is orthogonal to the central image plane, while the $z$-axis defines the depth of the volume, and the $x$-axis defines the lateral direction.

\subsubsection{Data}
\label{subsec:data}

\begin{figure}
    \centering
    \includegraphics[width=0.49\linewidth]{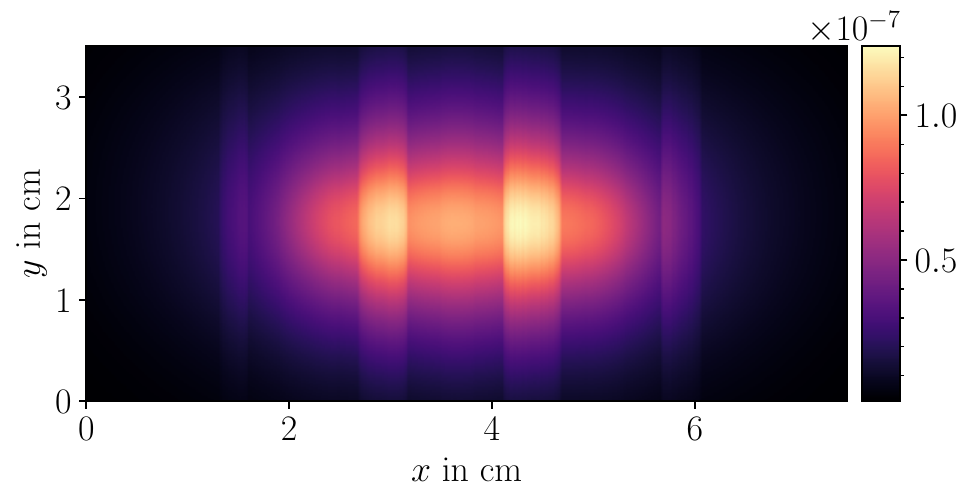}
    \includegraphics[width=0.5\linewidth]{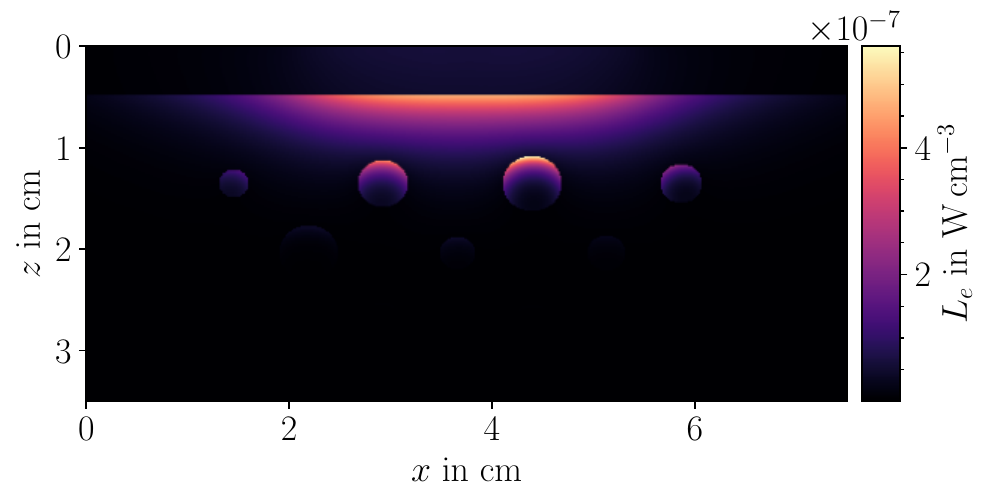}
    \includegraphics[width=0.49\linewidth]{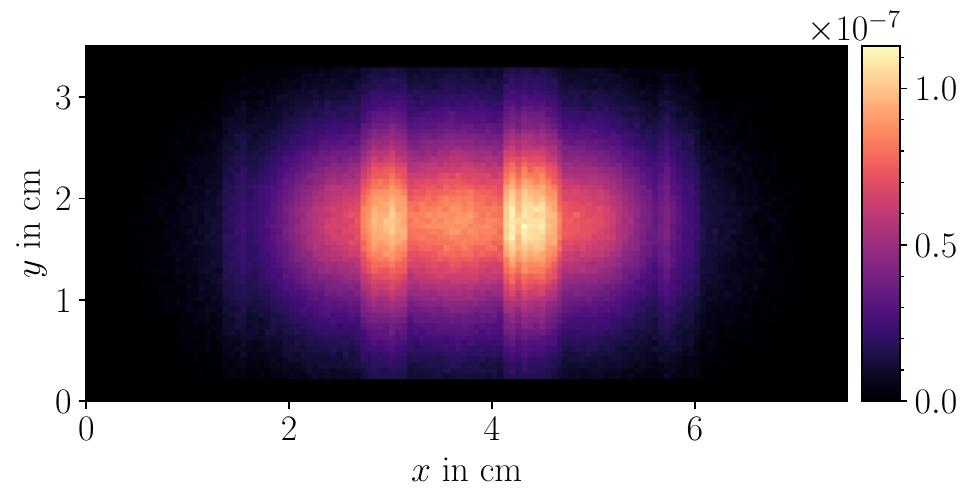}
    \includegraphics[width=0.5\linewidth]{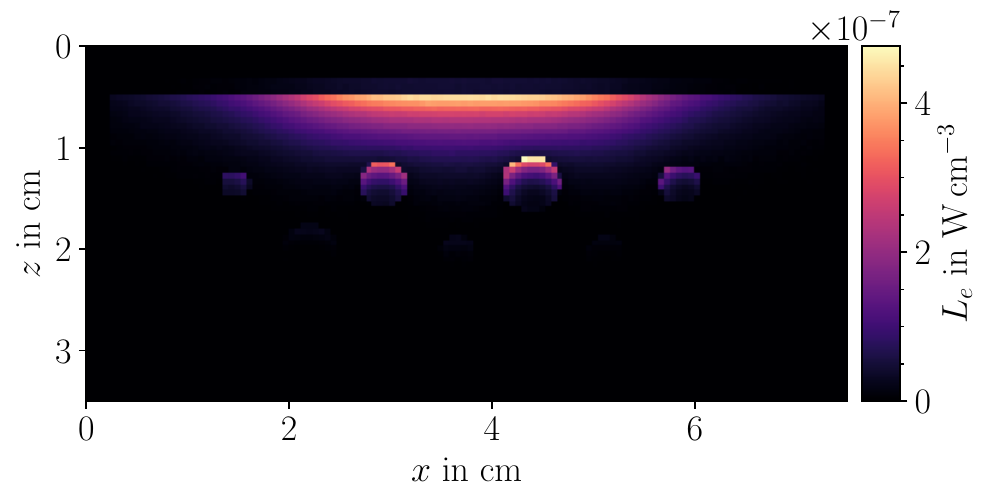}
    \caption{Depth (left) and orthogonal (right) projections of the absorbed energy density $L_e$ at full resolution (top) and as used in the reconstruction (bottom). The reconstruction dataset is masked to 0.2\,cm depth on all edges.}
    \label{fig:DataProj}
\end{figure}

The main dataset for this work was generated by solving the radiative transfer equation using a Monte Carlo method, yielding $L_e$ together with the optical coefficients. The spatial domain spans from $(0,0,0)$ to $(7.5,3.5,3.5)$\,cm in the $(x,y,z)$ directions with $(512,256,256)$ grid points, respectively. In Fig.~\ref{fig:DataProj} (top), we show depth and orthogonal projections of the full-resolution field. For reconstruction, $L_e$ and the coefficient fields are trilinearly interpolated onto the $(128,64,64)$ grid spanning the same volume, and Gaussian noise with a voxel-dependent standard deviation is added,
\begin{equation}
    \label{eq:noise}
    \sigma_i = \max\!\left(0.05\,L_{e,i},\,0.01\,\max_j L_{e,j}\right),
\end{equation}
where $L_{e,i}$ denotes the noise-free value at voxel $i$. The coarsened, noisy dataset used in the inference is shown in the bottom panels of Fig.~\ref{fig:DataProj}.

\begin{figure}
    \centering
    \includegraphics[width=0.7\linewidth]{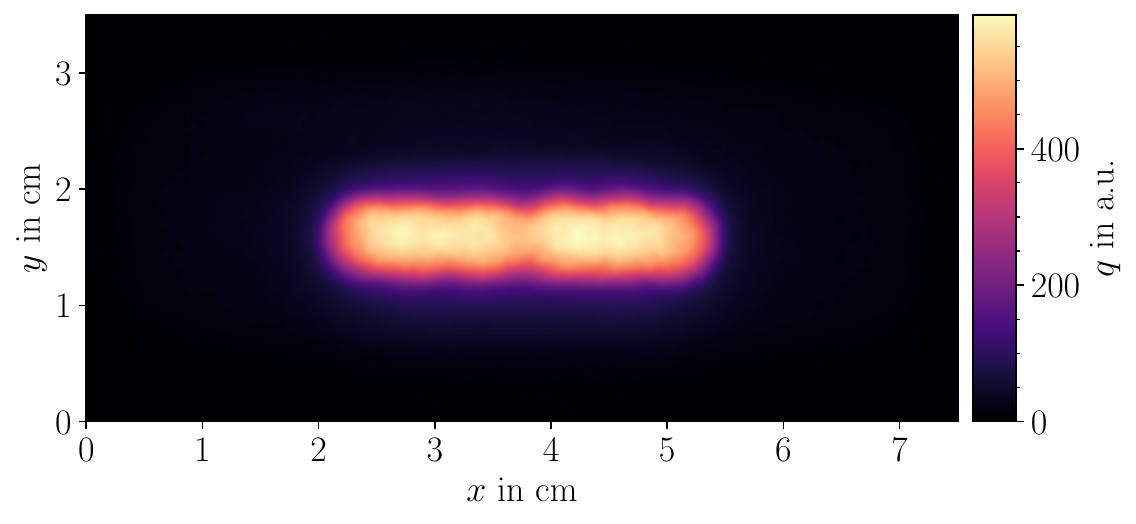}
    \caption{Illumination profile used to generate the dataset shown in Fig.~\ref{fig:DataProj}. The external source of the Monte Carlo setup is prescribed as a purely two-dimensional profile on the $z=0$ plane. We use this profile as the reference against which the structure of the reconstructed effective source is compared in Section~\ref{subsec:Recoresults}.}
    \label{fig:IllumProf}
\end{figure}

\begin{figure}
    \centering
    \includegraphics[width=0.51\linewidth]{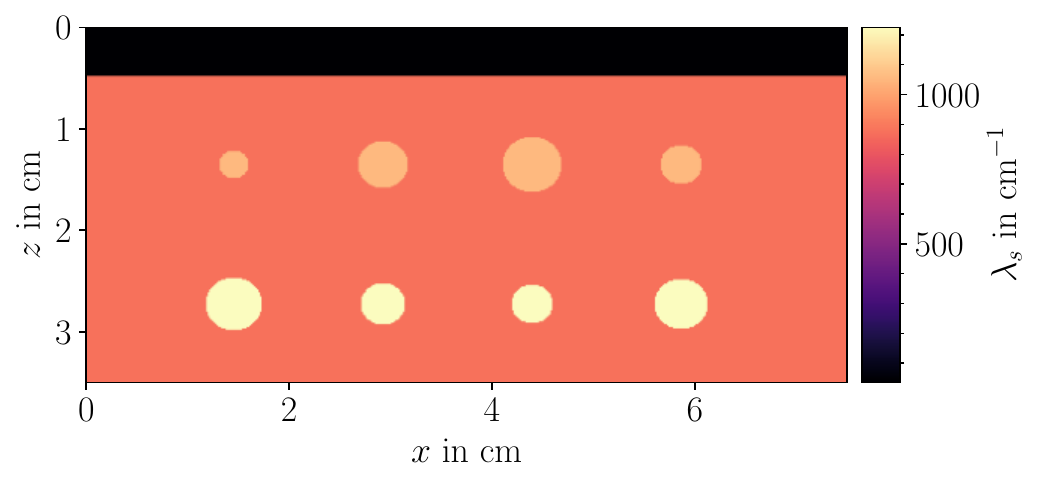}
    \includegraphics[width=0.48\linewidth]{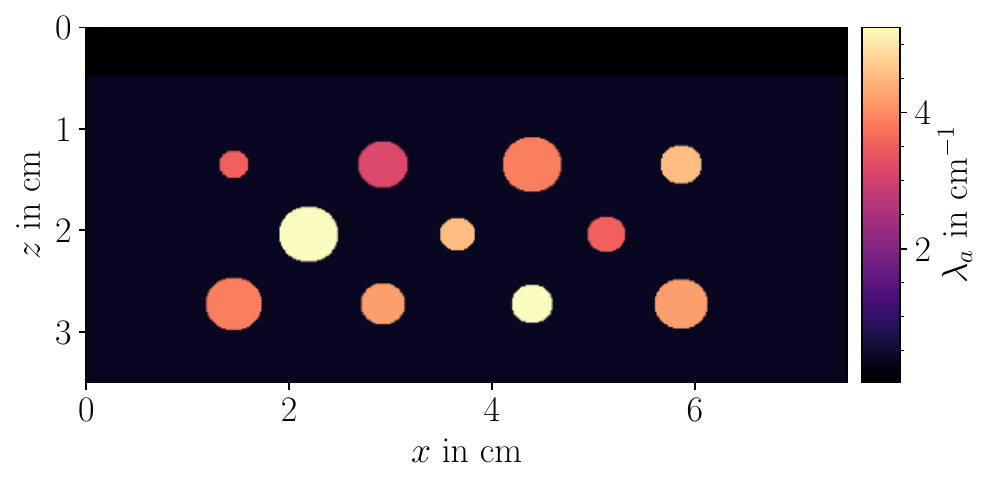}  
    \includegraphics[width=0.6\linewidth]{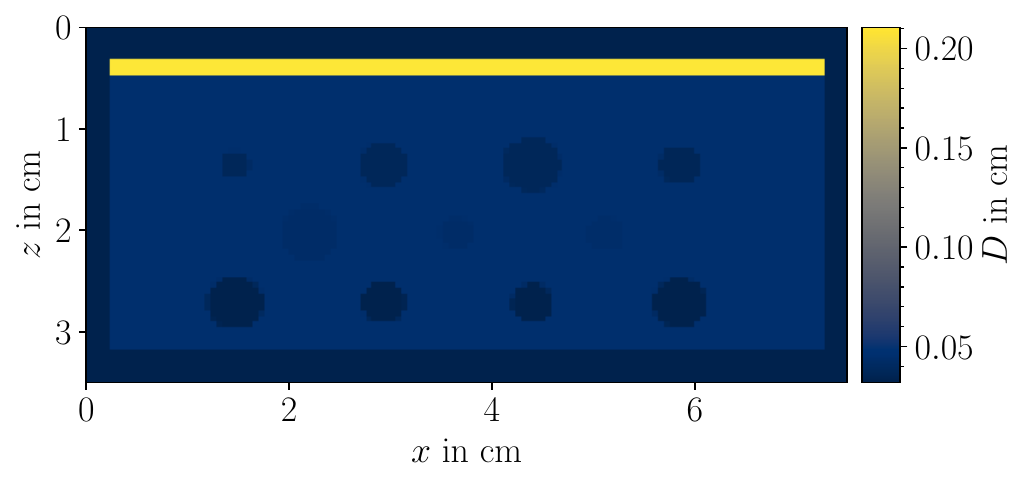}
    \caption{\textit{Top:} Orthogonal projections for $\lambda_s$ (left) and $\lambda_a$ (right). \textit{Bottom:} Orthogonal projection of the resulting diffusion coefficient $D$ as used in the reconstruction, also showing the masked region of depth 0.2\,cm at the edges.}
    \label{fig:medimcoeffs}
\end{figure}

In Fig.~\ref{fig:IllumProf}, we show the external light source used to generate this dataset. The internal tubular structures are introduced to mimic the presence of blood vessels within the tissue. These structures also locally influence the absorption and scattering coefficients. We show the orthogonal projections of $\lambda_a$ and $\lambda_s$ in Fig.~\ref{fig:medimcoeffs}, alongside the resulting diffusion coefficient $D$ (see Eq.~(\ref{eq:medimcoeffs})). The $L_e$ pattern in Fig.~\ref{fig:DataProj} includes two main regions: one external, located between the tissue and the radiation source; and one internal, which includes the blood-vessel-like structures. As part of the definition of the reduced scattering coefficient $\lambda_s'$, we define the scattering anisotropy $g$ as 0.9 internally and 0.5 externally. Since the radiation source is located at the top of the external region, at the origin of the $z$-axis in the volume over which the data are defined, we note that this region is dominated by the ballistic regime for radiative transport. Throughout this work, we therefore have to account for the error introduced by using the diffusion approximation in regions where the ballistic regime dominates (external region) and where Dirichlet boundary conditions modulate the resulting absorption distribution (volume edges). We account for these sources of error by masking the data up to a depth of 0.2\,cm on all edges of the grid volume.

\begin{figure}
    \centering
    \includegraphics[width=\linewidth]{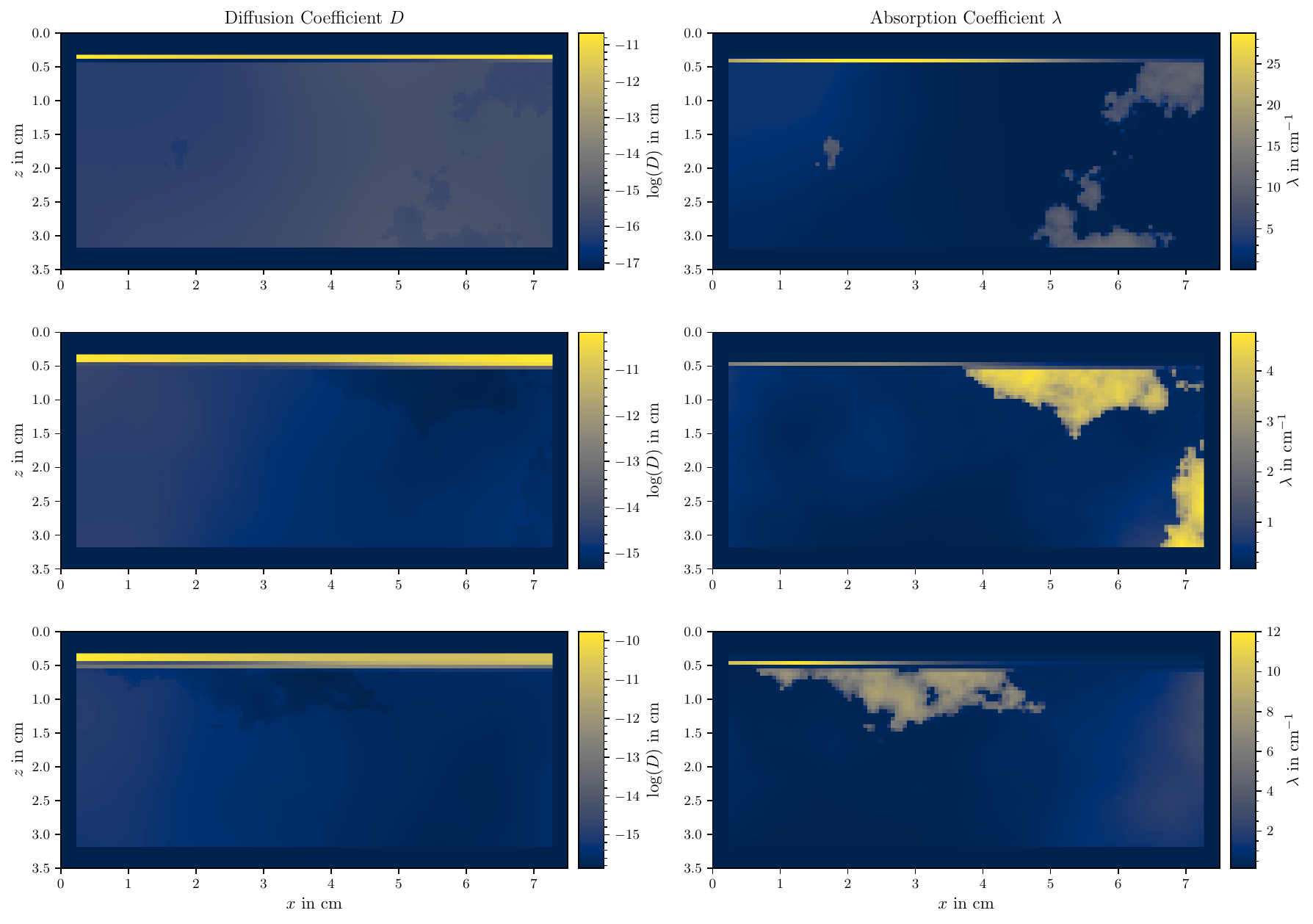}  
    \caption{Orthogonal projections of the diffusion coefficient $D$ (\textit{left}) and absorption coefficient $\lambda$ (\textit{right}) at full resolution for validation Datasets 17 (\textit{top}), 3 (\textit{middle}), and 28 (\textit{bottom}). Each row corresponds to one Gaussian random field realization of $\lambda_a$ and $\lambda_s$, and the ordering matches that of Fig.~\ref{fig:multi_validation}. All quantities are masked to 0.2\,cm depth on all edges.}
    \label{fig:valdata}
\end{figure}

To validate the capability of our code, we generate, in addition to the main dataset, 32 datasets with structures composed of Gaussian random fields through different realizations of the coefficients $\lambda_a$ and $\lambda_s$. The sets are generated with the same resolution and grid specifications as the reconstructed samples. In Fig.~\ref{fig:valdata}, we show the orthogonal projections of the resulting $\lambda$ and $D$ terms for three of the validation datasets, namely Datasets 17, 3, and 28. Each dataset contains a new absorbed energy density $L_e$ distribution, generated using the same external source as in Fig.~\ref{fig:IllumProf}. We compare these distributions against those obtained using the posterior samples of the inferred source distribution. Across all datasets, the same internal and external regions are defined for consistency, as the locations of the external source and the interface between media are fixed. However, the interface itself can vary in shape and location due to the relative values in $\lambda_a$ and $\lambda_s$. See Section~\ref{subsec:Recoresults} for more details on this validation process.

\subsubsection{Forward Model}
\label{subsec:forwardmodel}

For our reconstruction, the mapping defining the forward model has a simple structure: a source sample is used as input for DMGS (built using $D$, $\lambda$, and the grid specifications from the dataset). DMGS then provides a solution for the radiative fluence (i.e., $\rho$ in Eq.~(\ref{maineq})), which is multiplied by $\lambda$ to obtain $L_e$. Only voxels outside the masked edge region enter the comparison with the data, so we denote by $M$ the operator that selects these $N_{\mathrm{vox}}$ voxels and by $d$ the data vector collecting the corresponding noisy values of $L_e$ (Section~\ref{subsec:data}). The response of the model to a source $q$ is then

\begin{equation}
\label{eq:response}
    R(q) = M\left(\lambda\,S(q)\right),
\end{equation}
which maps a source realization to the space of the measured data.

\subsubsection{Prior Model}
\label{subsec:priormodel}

To regularize the reconstruction, we impose a structured prior on the source distribution. We express our prior knowledge about the effective source by setting a 3D profile. A prior realization of our source takes the form:

\begin{equation}
    q_{\rm{prior}} = c_1f_z(z)f_{xy}(x,y)\exp\left(c_2\phi(\vec{r})\right)
\end{equation}

Here, $c_1$ acts as a normalization constant; $f_{xy}(x,y)$ and $f_z(z)$ define the profile of the source distribution in the $xy$-plane and $z$-axis, respectively; and $c_2$ sets the scale of the Gaussian field $\phi(\vec{r})$.

For $f_{xy}(x,y)$, we define a 2D Gaussian profile centered in the $xy$-plane having the form:

\begin{equation}
    f_{xy}(x,y) = \exp\left(-0.5\cdot\left[(x-x_{\text{cen}})/\sigma_x\right]^2\right)\times\exp\left(-0.5\cdot\left[(y-y_{\text{cen}})/\sigma_y\right]^2\right)
\end{equation}
where $x_{\text{cen}}$ and $y_{\text{cen}}$ are the grid centers in the $xy$ plane, and $\sigma_x$ and $\sigma_y$ represent the width of the profile around $x_{\text{cen}}$ and $y_{\text{cen}}$, respectively. 

Along the depth direction, we impose a window function defined as the difference of two logistic functions:

\begin{equation}
    f_z(z) = \left[1+\exp(-2k\cdot\left(z-z_{\text{low}}\right))\right]^{-1}-\left[1+\exp(-2k\cdot\left(z-z_{\text{up}}\right))\right]^{-1}
\end{equation}

Here, $z_{\text{low}}=z_{\text{cen}}-0.5w$ and $z_{\text{up}}=z_{\text{cen}}+0.5w$ define the lower and upper edges, with $z_{\text{cen}}$ as the center and $w$ as the width of the window. Finally, $k$ represents the steepness at the edges of the window function, defining how sharp the decrease is outside of the window. 

\begin{figure}
    \centering
    \includegraphics[width=0.48\linewidth]{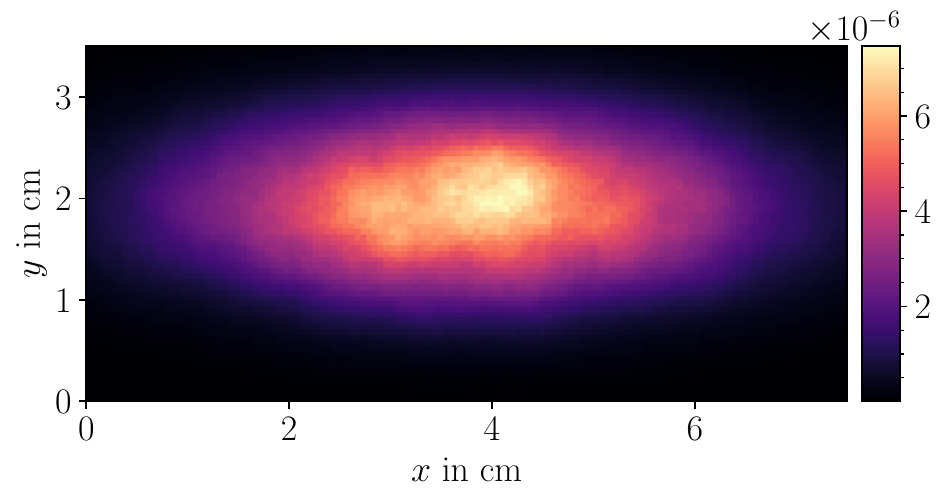}
    \includegraphics[width=0.5\linewidth]{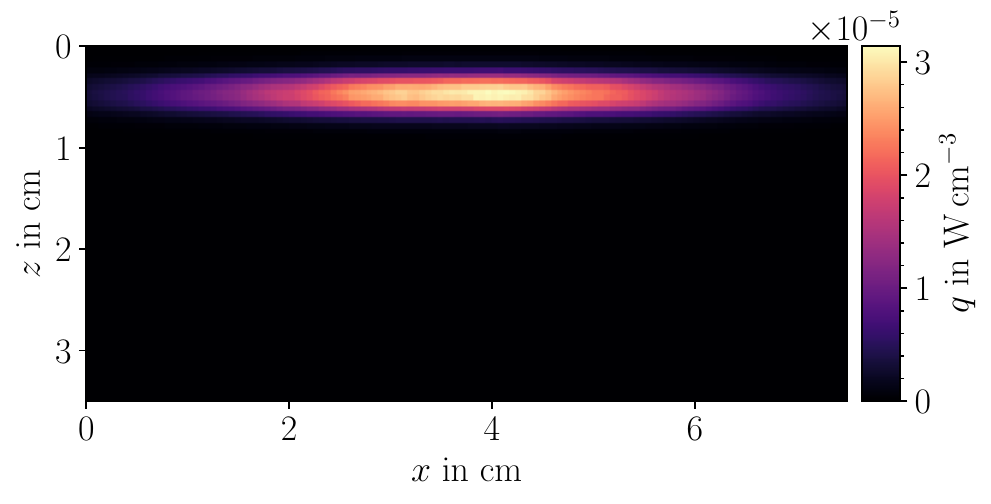}
    \caption{Orthogonal and depth projections of a source sample drawn from the prior model.}
    \label{fig:PriorSample}
\end{figure}

In Fig.~\ref{fig:PriorSample}, we show a prior sample that illustrates the actual 3D structure.

The free parameters of this prior fall into three groups, matching the columns of Table~\ref{tab:priorparams}: profile parameters ($c_1$, the centers and widths of $f_{xy}$ and $f_z$, and the edge steepness $k$), the field amplitude $c_2$, and the kernel parameters of $\phi(\vec{r})$. The field $\phi$ is a spatially correlated Gaussian field on the reconstruction grid with a Mat\'ern-like kernel specified by a spectral cutoff $\kappa_0$ and a neg-log-log slope $s_\phi$. All of these quantities are assigned log-normal priors in NIFTy.re and inferred jointly; Table~\ref{tab:priorparams} lists the prior means and standard deviations together with the corresponding posterior summaries.

\begin{table}[htbp]
\centering
\caption{Prior and posterior parameters of the source model, grouped into profile, field, and kernel parameters. Prior means and standard deviations specify log-normal priors in NIFTy.re; posterior summaries are computed from the stored VI samples reported in Section~\ref{subsec:Recoresults}. We discuss the interpretation of the posterior values in Section~\ref{subsec:Recoresults}.}
\label{tab:priorparams}
\begin{tabular}{llcccc}
\toprule
Group & Parameter & Prior mean & Prior std & Post.\ mean & Post.\ std \\
\midrule
\multirow{8}{*}{Profile}
 & $c_1$ & $2.5\times10^{-5}$ & $5.0\times10^{-6}$ & $2.56\times10^{-5}$ & $7.16\times10^{-6}$ \\
 & $x_{\mathrm{cen}}$ (cm) & 3.75 & 0.10 & 3.71 & 0.018 \\
 & $\sigma_x$ (cm) & 1.8 & 0.5 & 1.58 & 0.028 \\
 & $y_{\mathrm{cen}}$ (cm) & 1.75 & 0.10 & 1.75 & 0.065 \\
 & $\sigma_y$ (cm) & 0.8 & 0.1 & 1.43 & 0.096 \\
 & $w$ (cm) & 0.2 & 0.2 & 0.038 & 0.011 \\
 & $z_{\mathrm{cen}}$ (cm) & 0.5 & 0.1 & 0.17 & 0.023 \\
 & $k$ & 10 & 1 & 18.1 & 1.6 \\
\addlinespace
Field & $c_2$ & 0.10 & 0.01 & 0.085 & 0.0056 \\
\addlinespace
\multirow{2}{*}{Kernel}
 & $\kappa_0$ (cutoff) & 1.2 & 0.2 & 0.32 & 0.06 \\
 & $s_\phi$ (slope) & 3.5 & 1.0 & 8.8 & 1.7 \\
\bottomrule
\end{tabular}
\end{table}

\subsubsection{Likelihood and Inference Setup}
\label{subsec:inference}

Since the noise in Eq.~(\ref{eq:noise}) is uncorrelated between voxels, we collect its standard deviations in the diagonal covariance $N=\mathrm{diag}(\sigma_i^2)$ and take a Gaussian likelihood for the data given a source realization $q$,

\begin{equation}
\label{eq:likelihood}
    \mathcal{P}(d|q) \propto \exp\left[-\frac{1}{2}\sum_{i=1}^{N_{\mathrm{vox}}}\frac{\left(d_i-R_i(q)\right)^2}{\sigma_i^2}\right],
\end{equation}
with $R$ the response of Eq.~(\ref{eq:response}). As the noise level is known by construction, $N$ is fixed during the inference. Together with the prior of Section~\ref{subsec:priormodel}, Eq.~(\ref{eq:likelihood}) defines the posterior that we approximate with MGVI, where every sample and gradient evaluation calls DMGS in the fixed-iteration mode of Section~\ref{subsec:fulladjoint}.

We quantify the quality of the resulting fit by the reduced chi-squared statistic

\begin{equation}
\label{eq:chi2}
    \chi^2 = \frac{1}{N_{\mathrm{vox}}}\sum_{i=1}^{N_{\mathrm{vox}}}\frac{\left(d_i-R_i(q)\right)^2}{\sigma_i^2},
\end{equation}
which approaches unity when the residuals are consistent with the assumed noise level.

\section{Results}
\label{Results}

In the following, we demonstrate the practical compatibility of our differentiable multigrid solver with the push-forward operations required for gradient-based Bayesian inference. This is illustrated via the reconstruction of an effective diffuse radiative source using DMGS as part of the forward model. In addition, we report performance benchmarks comparing the hand-derived adjoint implementation (DMGS) with the autodiff-based implementation (JMGS).

\subsection{Reconstruction Quality}
\label{subsec:Recoresults}

Because the forward model is close to linear, the KL objective stabilizes after about one MGVI iteration. We report results from iteration 4 with eight stored posterior samples; continued resampling did not change these estimates.

\begin{figure}
    \centering
    \includegraphics[width=0.49\linewidth]{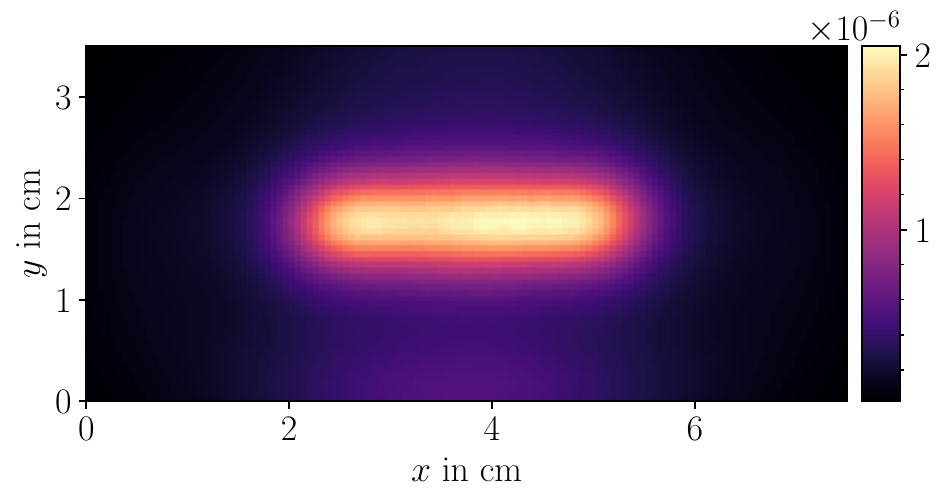}
    \includegraphics[width=0.5\linewidth]{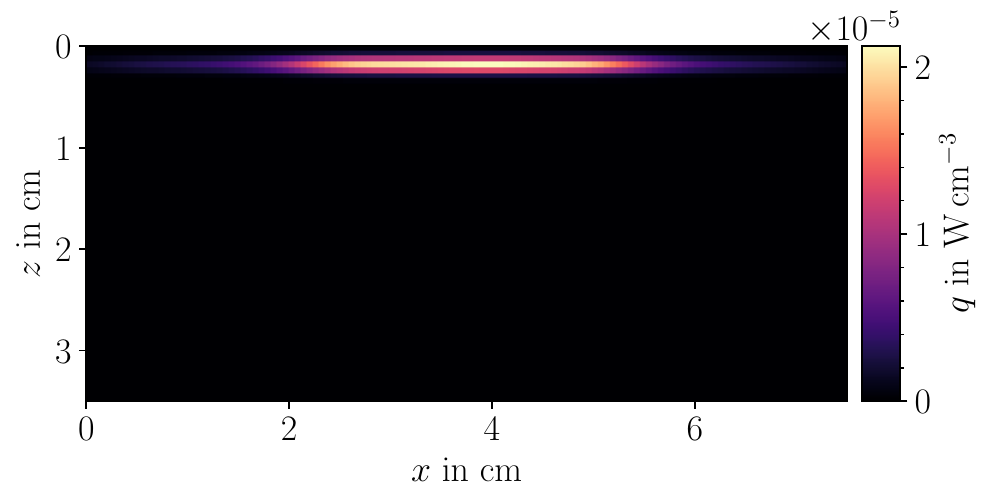}
    \caption{Depth (left) and orthogonal (right) projections of the reconstructed mean source.}
    \label{fig:Postmean}
\end{figure}

Fig.~\ref{fig:Postmean} shows orthogonal and depth projections of the posterior mean source distribution. The reconstructed source captures the dominant spatial structure of the illumination profile. In particular, the $xy$-plane morphology deviates from the initial Gaussian prior and instead reflects the structure of the true illumination (Fig.~\ref{fig:IllumProf}). Along the depth direction, the reconstruction places the effective source close to the imposed masking boundary, compensating for the excluded ballistic region.

Due to the diffusive nature of the forward model and the downscaling of the problem, the reconstruction primarily recovers large-scale features such as the location and extent of the source. Smaller-scale structures are smoothed out, as expected under the diffusion approximation. In addition, low-amplitude $L_e$ structures appear outside the main source region, particularly above about $2\,\mathrm{cm}$ and below about $1\,\mathrm{cm}$ in depth. These components compensate for residual ballistic contributions near the boundaries and allow the model to reproduce measurements in regions where the diffusion approximation is only partially valid.

\begin{figure}
    \centering
    \includegraphics[width=0.49\linewidth]{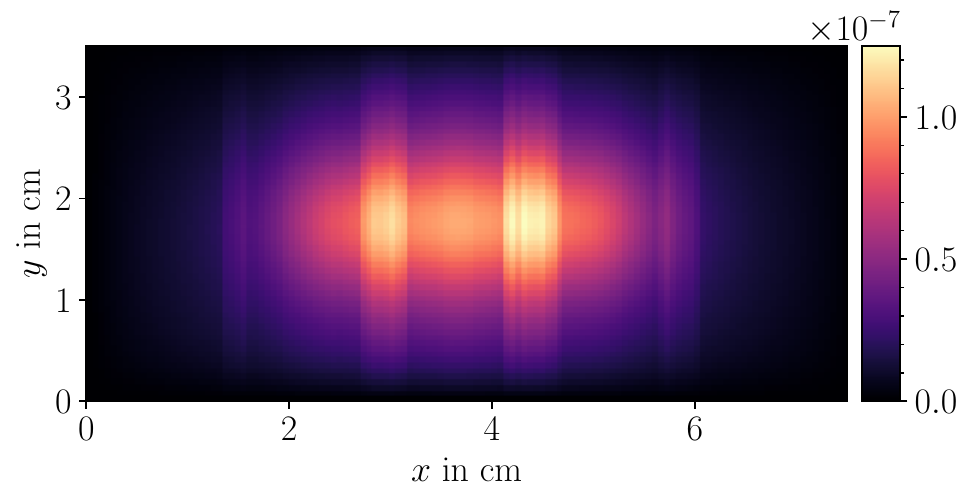}
    \includegraphics[width=0.5\linewidth]{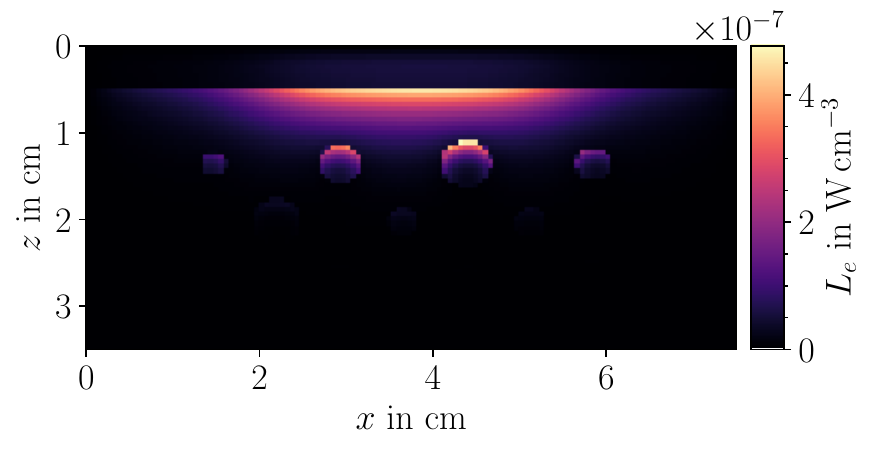}
    \caption{Depth (left) and orthogonal (right) projections for the response of the reconstructed mean source.}
    \label{fig:PostResponse}
\end{figure}

Fig.~\ref{fig:PostResponse} shows the forward response of the reconstructed mean source, i.e., the predicted absorbed energy density $L_e$. The reconstructed solution reproduces the observed data with high fidelity, yielding $\chi^2 \approx 1.1$ (see Eq.~(\ref{eq:chi2})). Deviations primarily occur near the boundaries of the reconstruction volume, where ballistic contributions are strongest.

\begin{figure}
    \centering
    \includegraphics[width=0.32\linewidth]{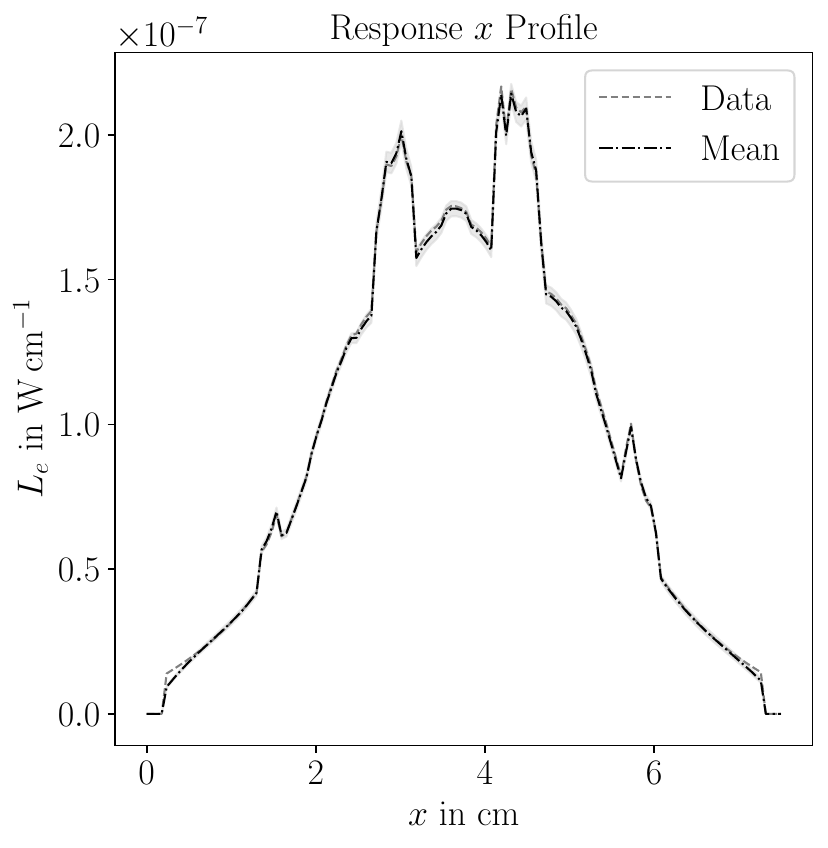}
    \includegraphics[width=0.305\linewidth]{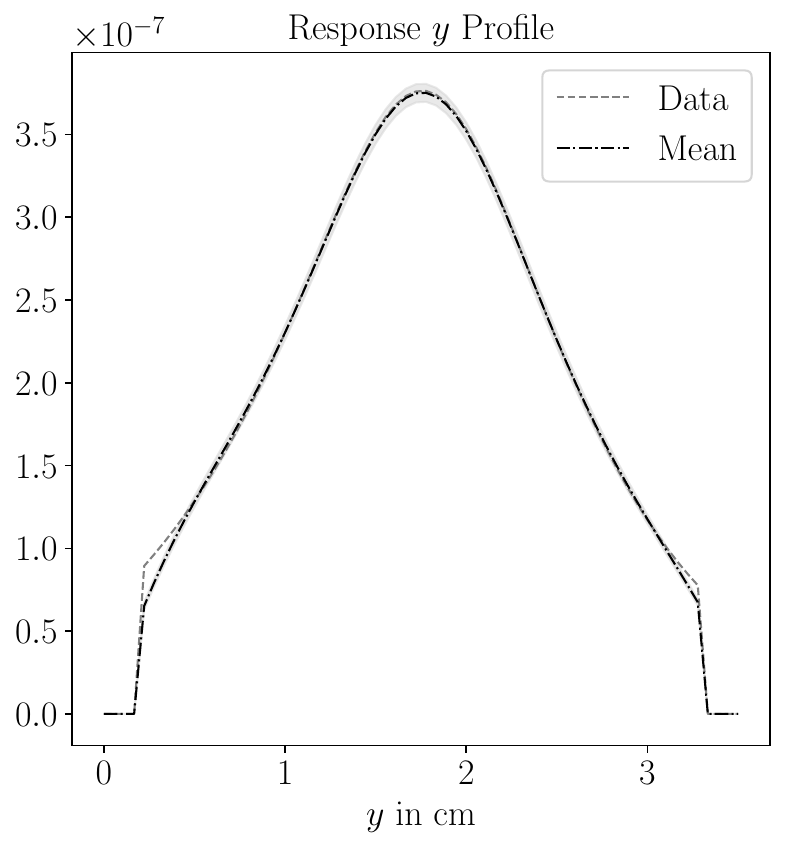}
    \includegraphics[width=0.31\linewidth]{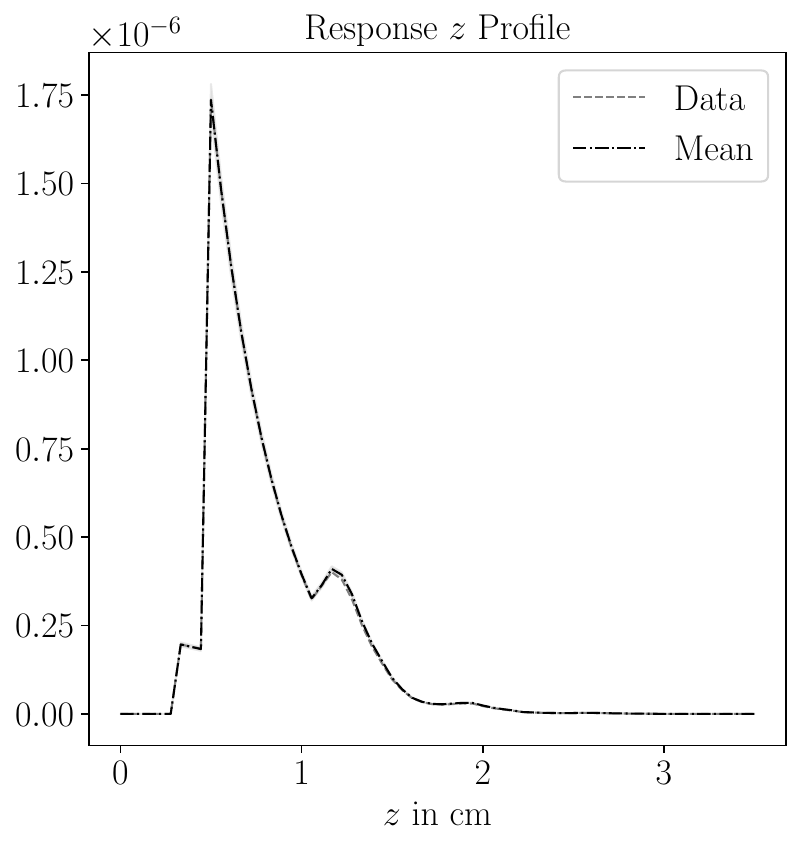}
    \caption{Profiles of the expected $L_e$ along the $x$-, $y$-, and $z$-axes for the posterior mean reconstruction compared to the data. The shaded region indicates the $1\sigma$ posterior uncertainty estimated from the sample dispersion.}
    \label{fig:Postprofile}
\end{figure}

A more detailed comparison is provided in Fig.~\ref{fig:Postprofile}, which shows profiles of the reconstructed response along each spatial axis. The reconstructed model accurately reproduces the measured data across all directions, with discrepancies again confined to regions near the boundaries.

The posterior parameters in Table~\ref{tab:priorparams} are consistent with this reconstruction picture. The lateral centers $x_{\mathrm{cen}}$ and $y_{\mathrm{cen}}$ remain close to their prior means, while their posterior standard deviations shrink, so that the source location is well determined. The lateral widths change more strongly: $\sigma_x$ contracts and $\sigma_y$ broadens to $1.43\,\mathrm{cm}$, yielding an anisotropic footprint that matches the $xy$-plane morphology of the true illumination profile. Along the depth direction, the window shifts toward the surface ($z_{\mathrm{cen}} = 0.17\,\mathrm{cm}$) rather than remaining at the prescribed tissue interface location, placing the effective source at the true depth of the illuminating source. The window also narrows ($w = 0.038\,\mathrm{cm}$) and sharpens ($k = 18.1$), so that the source sits against the masking boundary and can compensate for residual ballistic contributions.

The field and kernel parameters control the multiplicative correction to this profile. The field amplitude $c_2$ decreases slightly to $0.085$ with a tight posterior spread, so the profile carries most of the signal and the correction remains modest. For the Mat\'ern-like kernel, the spectral cutoff $\kappa_0$ drops from $1.2$ to $0.32$ and the neg-log-log slope $s_\phi$ steepens from $3.5$ to $8.8$. The posterior field is therefore smoother and more large-scale than assumed a priori, as expected for a diffusive forward model that suppresses fine spatial structure and leaves small-scale field fluctuations weakly constrained.

\subsection{Validation on Independent Data}

We evaluate generalization performance using 32 independent validation datasets. For each dataset, we compute the predicted $L_e$ using the known optical parameters ($D$, $\lambda$) alongside posterior samples of the reconstructed source. Values of $L_e$ are normalized to the total absorbed energy within the simulated volume. Elevated $L_e$ levels, on the order of $10^{-4}$ in normalized units, correlate with regions of high absorption. These levels drop rapidly as the distance from absorbing structures increases and the diffusion approximation becomes more reliable.

\begin{figure}[htbp]
\centering
\includegraphics[width=0.7\textwidth]{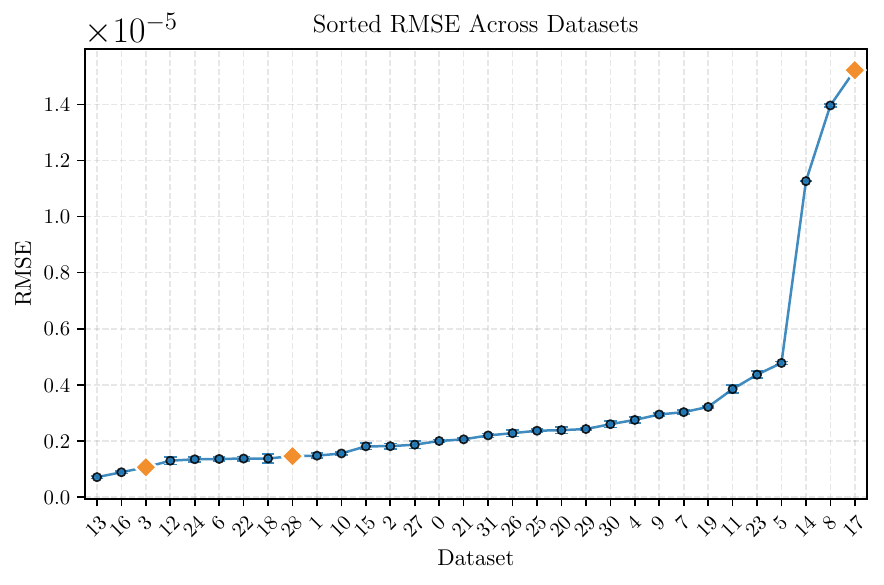}
\caption{Dataset-averaged root-mean-square error (RMSE) for all 32 validation datasets, sorted from lowest to highest RMSE (left to right) while retaining the original dataset indices on the horizontal axis. Datasets 3, 17, and 28 are marked with orange diamonds to contrast typical performance (Datasets~3 and 28) with the high-RMSE outlier (Dataset~17).}
\label{fig:rmse_errorbar}
\end{figure}

To quantify agreement, we calculate the voxel-wise root-mean-square error (RMSE) between predicted and reference $L_e$ on the masked reconstruction volume (same $0.2\,\mathrm{cm}$ edge mask as in the reconstruction dataset),
\begin{equation}
    \mathrm{RMSE}=\sqrt{\frac{1}{N_{\mathrm{vox}}}\sum_{i=1}^{N_{\mathrm{vox}}}\left(\hat{L}_{e,i}-L_{e,i}\right)^2},
\end{equation}
where $\hat{L}_{e,i}$ and $L_{e,i}$ are normalized values and $N_{\mathrm{vox}}$ is the number of unmasked voxels. We report dataset-averaged RMSEs in Fig.~\ref{fig:rmse_errorbar}. The majority of samples show excellent agreement on the order of $10^{-6}$, with three distinct outliers exhibiting higher errors on the order of $10^{-5}$, an order of magnitude above the typical case.

\begin{figure}[htbp]
\centering
\includegraphics[width=\textwidth]{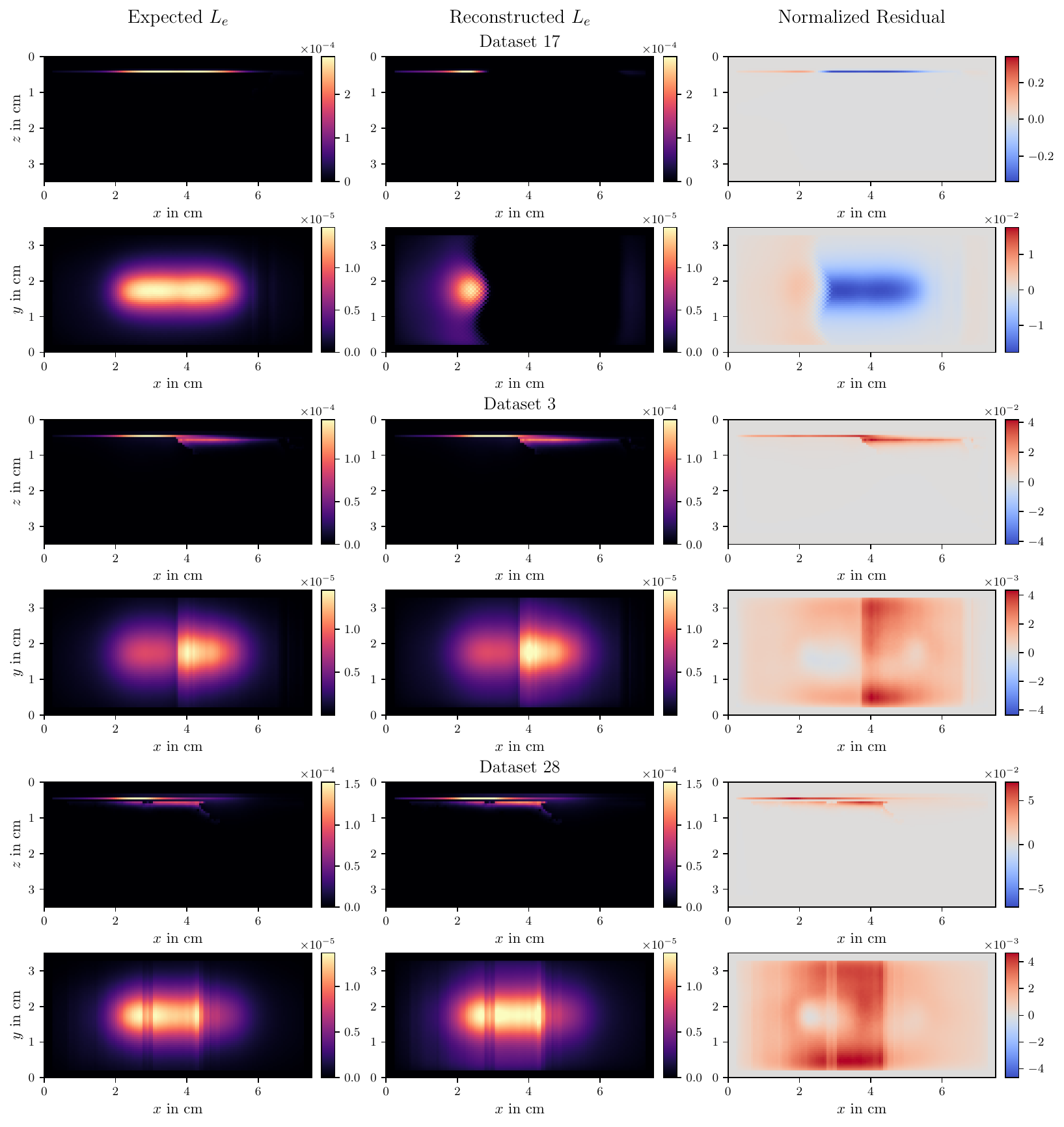}
\caption{Expected $L_e$ (left), reconstructed $L_e$ min--max rescaled to the reference range (middle), and residuals normalized by that reference range, $(\tilde{L}_e-L_e)/(L_{e,\max}-L_{e,\min}+\varepsilon)$ (right), for Datasets 17 (top), 3 (middle), and 28 (bottom). Displayed fields are clipped at zero to remove unphysical negative values in Dataset~17. The normalized residuals peak in Dataset~17, where high interface absorption makes the diffusion approximation unreliable.}
\label{fig:multi_validation}
\end{figure}

To further investigate these variations, we analyze three specific cases: Datasets 3, 17, and 28. Dataset~17 represents one of the high-RMSE outliers, while Datasets~3 and 28 are included for comparison to demonstrate typical, high-accuracy model performance. As illustrated by the $D$ and $\lambda$ projections in Fig.~\ref{fig:valdata}, the reconstructed source distribution effectively reproduces $L_e$ for different spatial combinations of diffusion and absorption, provided the diffusion approximation holds. Figure~\ref{fig:multi_validation} compares spatial distributions for these cases: reconstructed $L_e$ is min--max rescaled to the dynamic range of the Monte Carlo reference before plotting, and residuals are $(\tilde{L}_e-L_e)/(L_{e,\max}-L_{e,\min}+\varepsilon)$, with $\tilde{L}_e$ the rescaled prediction and $L_{e,\min}$, $L_{e,\max}$ the reference extrema over the masked volume.

However, the forward-model prediction degrades when the diffusion approximation becomes unreliable. This is visible in Dataset~17, which features a comparatively high absorption coefficient concentrated near the interface of the internal and external regions. There, deviations reach upwards of 0.2 in the range-normalized residuals. In this case, residual ballistic contributions dominate the highly absorbing interface; the mismatch reflects a limitation of the diffusion forward model rather than of the differentiable solver or adjoint implementation. Conversely, configurations like Datasets~3 and 28, where absorption is localized deeper within the medium, yield much better agreement and strongly align with the assumptions of our forward model.

We now shift our focus from evaluating statistical accuracy to assessing computational efficiency. The following section details the performance benchmarks of the actual code.

\subsection{Performance Benchmarks}

We compare the computational performance of the DMGS and JMGS implementations in two dimensions and report the performance of DMGS in three dimensions. Benchmarks are performed using identical solver configurations, differing only in the implementation of the forward and adjoint operations. All timings use JAX~0.9.0 in float64 on dual 32-core Intel Xeon~8358 (Ice Lake) CPUs; reported runtimes exclude one-time JIT compilation and memory denotes peak resident set size (RSS) per call.

In all cases, we evaluate three modes: direct solver execution ("Solve"), forward-mode differentiation (JVP), and reverse-mode differentiation (VJP). Tests are conducted over grid resolutions from $2^{19}$ to $2^{26}$ in 2D and up to $2^{28}$ in 3D. Each configuration is repeated 20 times to account for runtime variability.

\begin{figure}
    \centering
    \includegraphics[width=0.45\linewidth]{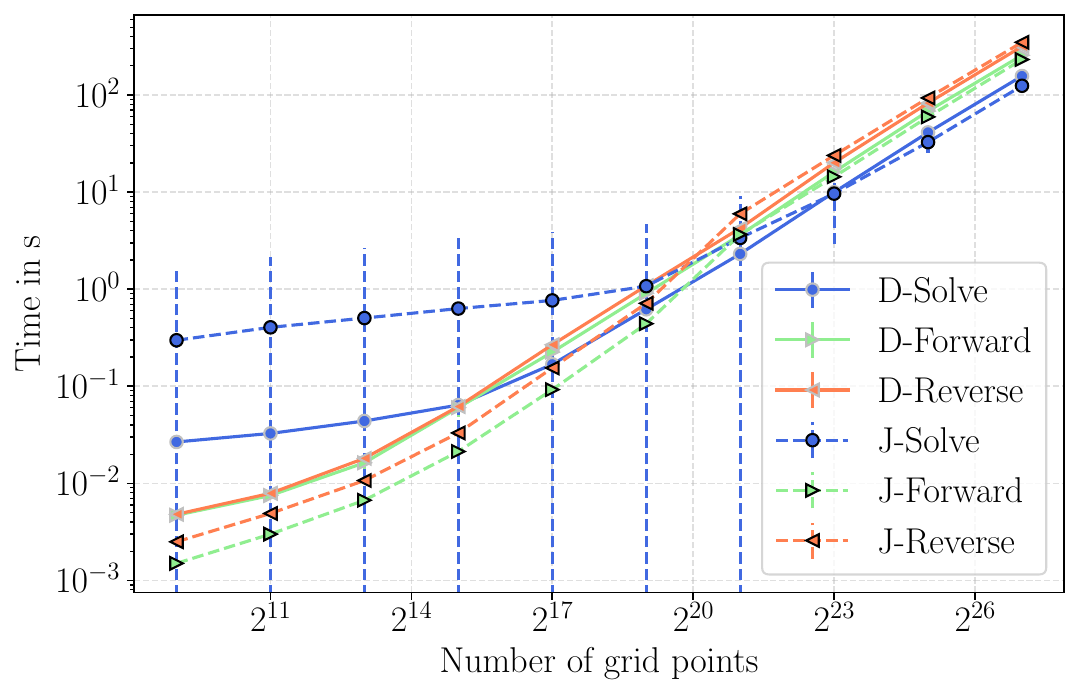}
    \includegraphics[width=0.45\linewidth]{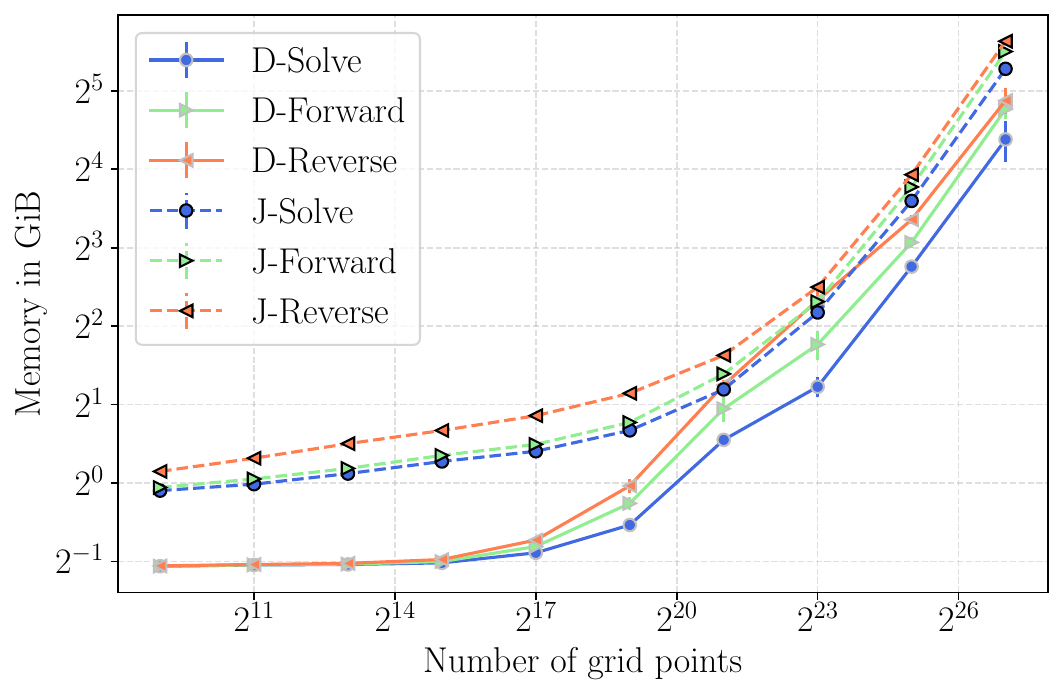}
    \includegraphics[width=0.45\linewidth]{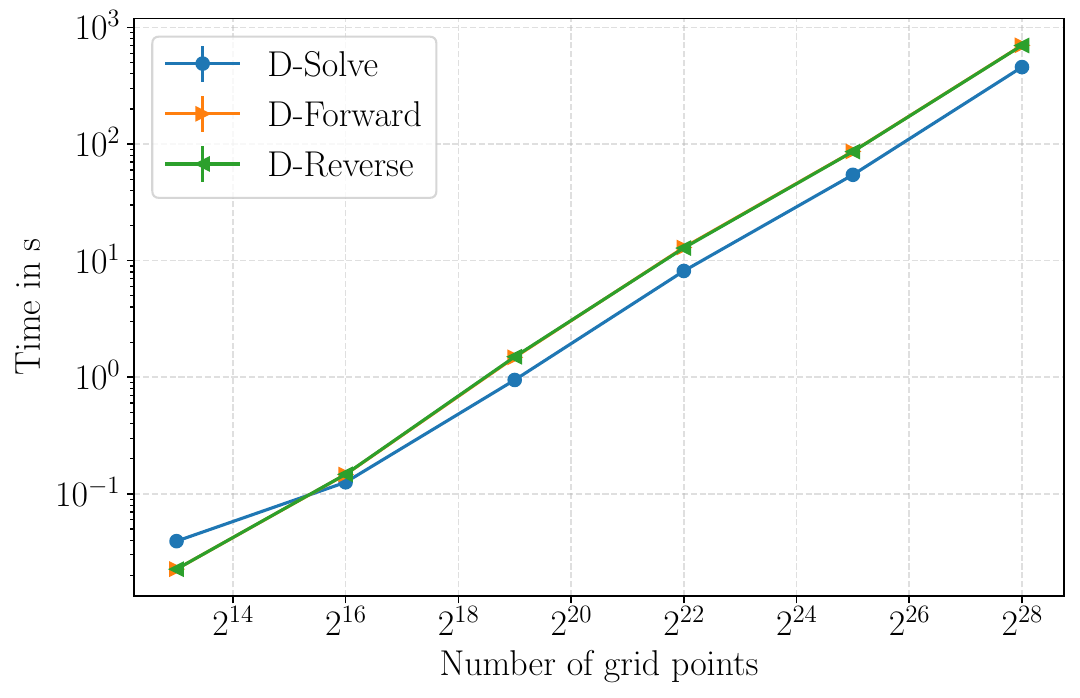}
    \includegraphics[width=0.45\linewidth]{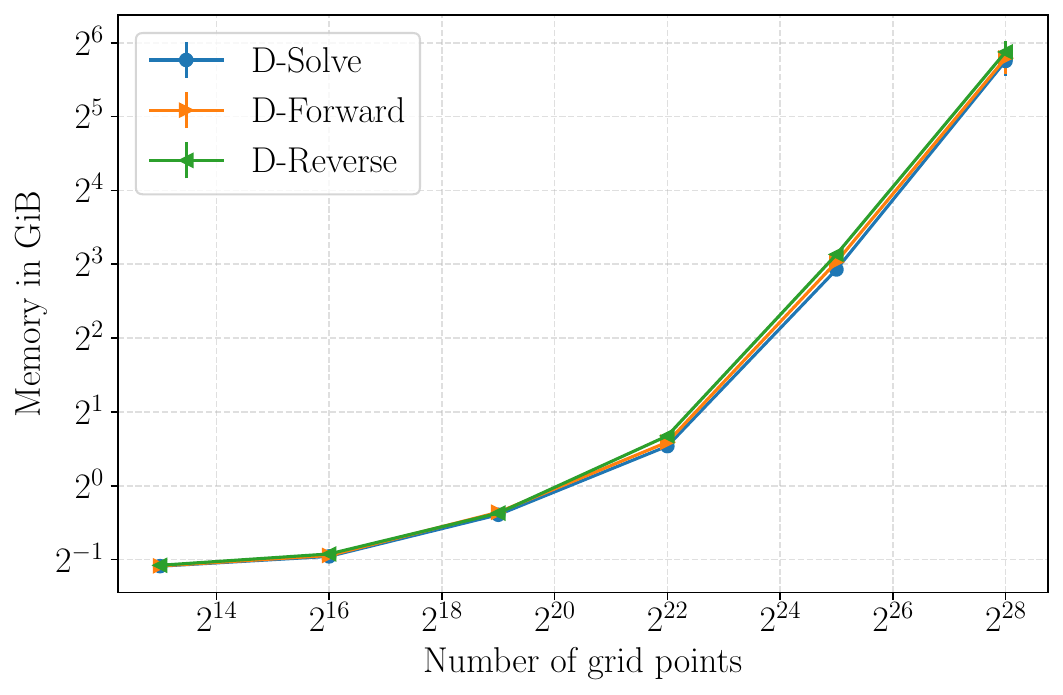}
    \caption{Benchmark runtime (left) and peak memory (right) versus grid resolution in 2D (top) and 3D (bottom). Main result: DMGS matches JMGS in runtime but uses less memory in reverse mode by avoiding storage of the full AD trace. Error bars show the standard deviation over 20 runs (JAX~0.9.0, dual Intel Xeon~8358 CPUs, float64; JIT excluded).}
    \label{fig:Benchmarks}
\end{figure}

Fig.~\ref{fig:Benchmarks} summarizes runtime and memory as a function of grid resolution. Both implementations exhibit comparable runtime at high resolution; the main difference is memory, where DMGS avoids storing the full JAX trace in reverse mode.

The primary difference between the implementations lies in memory consumption. While both approaches exhibit similar runtime scaling, the autodiff-based JMGS implementation requires storing intermediate states during reverse-mode differentiation, resulting in consistently higher memory usage. In contrast, DMGS avoids storing the full computational trace and therefore achieves lower memory consumption across all tested configurations. In 2D, reverse mode is up to $24\%$ slower than forward mode for DMGS and up to $67\%$ slower for JMGS; in 3D, DMGS forward and reverse modes differ by less than 1\% in runtime.

These results highlight the trade-off between flexibility and efficiency: autodiff-based implementations provide a convenient framework for constructing differentiable solvers, whereas hand-derived adjoint implementations offer significant advantages in memory efficiency, particularly for large-scale three-dimensional problems.

\section{Conclusions}
\label{Conclusions}

We have developed a differentiable multigrid solver (DMGS), implemented in C++ with explicitly derived adjoint operations, and demonstrated its applicability within gradient-based Bayesian inference. By integrating DMGS into the NIFTy framework, we establish that multigrid solvers can serve as efficient and fully compatible components of the push-forward mapping required for variational inference in PDE-constrained problems.

Applied to a medical imaging problem based on radiative transfer in organic tissue, the method successfully reconstructs an effective diffuse light source from data generated via the radiative transfer equation. The reconstructed source achieves a statistically consistent fit to the main dataset, with $\chi^2 \approx 1.1$. It accurately recovers the dominant spatial structure of the external illumination while compensating for systematic residuals near the boundaries. These residuals arise from the diffusion approximation, which shifts the reconstructed source toward the top layers of the domain to account for unmodeled ballistic photon transport contributions.

The reconstruction generalizes well to independent data. Across 32 validation datasets with varying optical properties, the predicted $L_e$ distributions remain consistent with Monte Carlo simulations. On average, the RMSE is an order of magnitude smaller than the local $L_e$. Larger deviations are confined to regions near domain boundaries and cases with strongly absorbing interfaces, where the diffusion approximation breaks down due to increased ballistic contributions.

We further compared the performance of DMGS with a JAX-based autodiff implementation (JMGS) using identical solver configurations. In two dimensions, reverse-mode differentiation is up to $24\%$ slower than forward mode for DMGS and up to $67\%$ slower for JMGS, with these differences decreasing at higher resolutions. Memory consumption shows a clearer distinction: reverse mode requires approximately $12\%$ more memory than forward mode for DMGS and $20\%$ more for JMGS. In three dimensions, DMGS exhibits minimal overhead, with average memory increases of $3\%$ and negligible runtime differences between forward and reverse modes.

These results demonstrate that while autodiff-based implementations provide flexibility, hand-derived adjoint solvers offer substantial advantages in memory efficiency without compromising runtime performance. This becomes increasingly relevant for highly resolved three-dimensional problems, where memory constraints become the limiting factor.

To our knowledge, this work is among the first integrations of a multigrid solver with explicit adjoint support into a metric-based variational inference framework. Future work includes nonlinear PDE solvers, differentiation with respect to $D$ and $\lambda$, larger GPU-accelerated 3D runs, and application to experimental optoacoustic data.

\begin{acknowledgements}
  All figures in this publication have been created using \emph{matplotlib}~\citep{Hunter2007}.
  This research was funded in part by the Austrian Science Fund (FWF) [I 5925-N].
  Part of this work was supported by the German \emph{Deut\-sche For\-schungs\-ge\-mein\-schaft, DFG\/} project number 495252601.
  Philipp Frank acknowledges funding through the German Federal Ministry of Education and Research for the project ErUM-IFT: Informationsfeldtheorie f\"ur Experimente an Gro{\ss}forschungsanlagen (F\"orderkennzeichen: 05D23EO1). Vo Hong Minh Phan acknowledges support from the Initiative Physique des Infinis (IPI), a research training program of the Idex SUPER at Sorbonne Universit\'e. Dominik J\"ustel and Philipp Haim acknowledge funding from the European Research Council (ERC) under the European Union's Horizon Europe research and innovation programme under grant agreement No.~101041936 (EchoLux), and from the Bavarian Ministry of Economic Affairs, Energy and Technology (StMWi) (DIE-2106-0005// DIE0161/02, DeepOpus).
\end{acknowledgements}

\bibliographystyle{spmpsci}
\bibliography{references}
\end{document}